\documentclass[reqno,a4paper,11pt]{amsart}
\usepackage{amsmath,amsthm,amssymb,mathrsfs}
\usepackage[
left = 2.25cm,
right = 2.25cm,
top = 2.2 cm,
bottom = 2.2 cm
]{geometry}
\usepackage{graphicx,xcolor,tikz,mathtools}
\usepackage{bm}
\usepackage{amsaddr}
\usepackage{enumitem}
\usepackage{esint} 
\usepackage{hyperref}

\usepackage{todonotes}
\usepackage{cite}

\makeatletter
\define@key{todonotes}{Helmut}[]{%
	\setkeys{todonotes}{author=\textbf{Helmut},inline,color=red!10}}%
\define@key{todonotes}{Andrea}[]{%
	\setkeys{todonotes}{author=\textbf{Andrea},inline,color=green!20}}%
\define@key{todonotes}{Yadong}[]{%
	\setkeys{todonotes}{author=\textbf{Yadong},inline,color=cyan!20}}%
\makeatother

\theoremstyle{plain}
\newtheorem{theorem}{Theorem}[section]

\newtheorem{definition}[theorem]{Definition}
\newtheorem{lemma}[theorem]{Lemma}
\newtheorem{proposition}[theorem]{Proposition}
\newtheorem{remark}[theorem]{Remark}
\theoremstyle{definition}

\numberwithin{equation}{section}

\newcommand{\bu}{\mathbf{u}}
\newcommand{\bv}{\vphi}
\newcommand{\bw}{\mathbf{w}}
\newcommand{\bn}{\mathbf{n}}

\renewcommand{\d}{\mathrm{d}}
\newcommand{\dx}{\,\d x}
\newcommand{\dt}{\,\d t}

\newcommand{\dtau}{\,\d \tau}

\newcommand{\norm}[1]{\left\Vert #1 \right \Vert}
\newcommand{\normm}[1]{\Vert #1 \Vert}
\newcommand{\normmm}[1]{\left\vert #1 \right\vert}

\def\Div{\mathrm{div}\,}

\newcommand{\dist}{\mathrm{dist}}

\allowdisplaybreaks[4]

\newcommand{\R}{\mathbb{R}}

\renewcommand{\d}{\mathrm{d}}

\def\Div{\mathrm{div}\,}

\def\Ld{L^2(\Omega)}
\def\bLd{\mathbf L^2(\Omega)}
\def\bLds{\mathbf L^2_\sigma(\Omega)}
\def\N{\mathbb N}
\def\R{\mathbb R}
\def\vphi{\varphi}
\def\Hu{H^1(\Omega)}

\begin{document}
	
	\title[Convergence to equilibrium of weak solutions to the Cahn--Hilliard equation]{Convergence to equilibrium of weak solutions\\ to the Cahn--Hilliard equation \\with non-degenerate mobility\\ and singular potential}
	
\author[M. Grasselli and A. Poiatti]{
	\small
Maurizio Grasselli$^\ast$ and
	Andrea Poiatti$^\ddagger$
}

\address{
	$^\ast$Dipartimento di Matematica, Politecnico di Milano, 20133 Milano, Italy
}
\email{maurizio.grasselli@polimi.it}

\address{
	$^\ddagger$Faculty of Mathematics, University of Vienna,
	1090 Vienna, Austria
}
\email{andrea.poiatti@univie.ac.at}

	
	

	\begin{abstract}
    We consider the classical initial and boundary value problem for the Cahn--Hilliard equation with non-degenerate mobility and singular (e.g., logarithmic) potential.
    We prove that any weak solution converges to a single equilibrium using \textit{only} minimal assumptions, that is, the existence of a global weak solution which satisfies an energy inequality. This result appears to be new in the literature and also holds in the three-dimensional case, which was an open problem due to the lack of regularity results, especially when the mobility is just a continuous function. We then prove the same result for a Cahn--Hilliard-Navier--Stokes type system with unmatched densities and viscosities proposed by Abels, Garcke, and Gr\"{u}n (Math.
Models Methods Appl. Sci. 22, 2012), always assuming a non-degenerate mobility. We expect that this novel method can be used to analyze the same issue for other models  where the regularization properties of the solutions are unknown or unlikely.
    \end{abstract}
	
	\maketitle
	
	
	\section{Introduction}
    The celebrated Cahn--Hilliard equation was proposed in \cite{CH1} (see also \cite{Cahn,CH2}) as a model to describe phase separation in binary alloys. Since then, it has been used in many different contexts which have in common a phase separation phenomenon (see, for instance, \cite{Miranville} and references therein). In particular, we recall that phase separation has recently become a paradigm in Cell Biology (see \cite{Alberti,Dolgin} and references therein). Therefore, we can say that the importance of the Cahn--Hilliard equation for applications has increased greatly during the last decades, as testified by the scientific literature. Referring to the original setting, assume that the alloy occupies a bounded domain $\Omega \subset \mathbb{R}^d$, $d\in\{1,2,3\}$, and $\vphi$ denotes the relative concentration of the difference of the two components. The Cahn--Hilliard equation can be written as follows (see, for instance, \cite{Elliott})
     \begin{align}
     \label{A1} &\partial_t\vphi-\Div(m(\vphi)\nabla \mu)=0,\quad\text{ in }\Omega\times(0,\infty)\\&
      \label{A2}  \mu=-\Delta\vphi+f'(\vphi),\quad\text{ in }\Omega\times(0,\infty),
    \end{align}
    where we have set some constants equal to the unity.
    Here, $m(\cdot)$ is the mobility, $\mu$ is the chemical potential (i.e., the first variation of the free energy with respect to $\varphi$), and $f$ is the potential density whose typical expression is the following (a.k.a. Flory--Huggins potential)
\begin{equation}
\label{FHP}
    f(s)=\frac{\theta}{2}((1+s)\text{ln}(1+s)+(1-s)\text{ln}(1-s)) -\frac{\theta_0}{2}s^2,
    \quad\forall\,s\in(-1,1).
\end{equation}
    We recall that $\theta>0$ and $\theta_0>0$ stand for the absolute temperature and a critical temperature (depending on the mixture), respectively. Also, we point out that phase separation takes place only if $\theta$ is below the critical temperature, that is, when $f$ is a double well.

    The standard initial and boundary value problem for equation \eqref{A1} is characterized by the following conditions
 \begin{align}
        &\partial_n\vphi=0,\quad\text{ on }\partial\Omega\times(0,\infty) \label{Af0}\\&
        \partial_n\mu=0,\quad\text{ on }\partial\Omega\times(0,\infty)\label{Af1}\\&
        \vphi(0)=\vphi_0,\quad\text{ in }\Omega.\label{Af}
    \end{align}
     These conditions ensure the conservation of the total mass, that is, $\int_\Omega \vphi(x,t)\dx$ is constant
     for all times $t\geq 0$. We recall that other relevant boundary conditions are the dynamic ones (see, for instance, the review \cite{Wu} on the topic).

    On account of its importance in a number of applications, problem \eqref{A1}-\eqref{Af} has been extensively studied from the numerical viewpoint.
    The theoretical results are many, but they depend on the nature of mobility as well as on the choice of $f$. Most of them are related to the constant mobility case where $f$ is often approximated by a fourth-order double well polynomial (a.k.a. smooth potential). However, in this case, one cannot ensure that $\vphi$ takes its values in the physical range $[-1,1]$.
    
    The main theoretical aspects (i.e., well-posedness, regularity, and longterm behavior) of the constant mobility case have been explored in detail. On the other hand, a significant issue is still open, namely, the validity of the instantaneous strict separation property in dimension three in the case of a logarithmic potential $f$ as in \eqref{FHP} (see, for instance, \cite{GalP} and references therein). The strict separation property means that $\vphi$ stays uniformly and instantaneously away from the pure phases $\pm 1$ so that $f^\prime$ is no longer singular.

    If mobility degenerates in the pure phases (see \cite{CH2}) and $f$ is the Flory--Huggins potential, then the only theoretical result is an existence theorem (see \cite{EG}). If mobility is non-degenerate, the picture is slightly better. Namely, there exists a global weak solution which satisfies an energy identity (see \cite{BB} and also \cite{Schimperna}). In dimension two, uniqueness, regularity results, and the validity of the strict separation property also hold (see \cite{BB,CGGG}).

    The convergence to equilibrium of solutions to \eqref{A1}-\eqref{Af} for constant mobility and smooth potential was first studied in \cite{RybHof}. Then, a great step forward was made in \cite{AW} for the case of singular potential (e.g. the Flory--Huggins one). This result was made possible by the asymptotic validity of the strict separation property. This fact is essential to ensure the real analyticity of the singular potential, an essential requirement in order to use a suitable version of the {\L}ojasiewicz--Simon inequality. We recall that Cahn--Hilliard (or Allen--Cahn) equations have a continuum of stationary states even in dimension two, so that standard techniques based on the Lyapunov functional cannot be applied  (see \cite[Remark 2.3.13]{Haraux}). Moreover, if the nonlinearity is not a real analytic function, then the solution might not converge to equilibrium (see \cite{PolSim} and references therein).

    In the case of non-degenerate mobility, convergence to equilibrium has recently been established only in dimension two in \cite{CGGG}, using the regularization properties of the solution, the energy identity, and the instantaneous strict separation. On the other hand, no result is known in dimension three because of the lack of regularity results. Here we propose a novel approach that only needs the existence of a global weak solution and the validity of an energy inequality. In addition, the mobility coefficient is assumed to be continuous only. Therefore, any weak solution to \eqref{A1}-\eqref{Af} converges to a single stationary state even in dimension three with minimal assumptions on mobility (just continuous in $[-1,1]$ and bounded from below by a positive constant). 
    
    However, there is more. Indeed, the fact that only two ingredients (i.e., existence of a global weak solution and validity of the energy identity) are needed opens the way to prove analogous results for more complicated models characterized by non-degenerate mobility. For instance, multi-component Cahn--Hilliard or Allen--Cahn equations (cf. \cite{GGPS,GPCAC} and their references), Cahn--Hilliard equations on surfaces evolving with prescribed velocity (see \cite{CEGP}, see also \cite{AGGP1,AGGP2} for some coupled problems), Cahn--Hilliard equation with dynamic boundary conditions (see \cite{Wu} and its references), and Cahn--Hilliard--Darcy systems (cf. \cite{Gio} and references therein). Recall that in the case of Cahn--Hilliard--Darcy, the issue is open even for constant mobility in dimension three. In conclusion, we point out that, even in the case of constant mobility, the present method noticeably simplifies the existing proofs. 
    
    In order to support the feasibility of our claim, as well as to show the robustness of our novel method, in addition to problem \eqref{A1}-\eqref{Af}, we will show an application to the so-called  Abels-Garcke-Gr\"{u}n (a.k.a. AGG) model with non-degenerate mobility, that is, a Cahn--Hilliard-Navier--Stokes system with unmatched densities and viscosities first introduced in  \cite{AGG}. This application also shows that other similar existing results for Cahn--Hilliard--Navier--Stokes systems, even in the case of constant mobility, can be significantly simplified (see, for instance, \cite{AGGG,AGGmulti}). Moreover, extensions to more complex systems (see, e.g., \cite{DiPG,HeWu})
    are also possible.
      
    The problem analyzed in \cite{ADG} is the following
    \begin{align}
       \label{A1AGG} & \partial_t(\rho(\vphi)\bu)+\Div(\bu\otimes \rho(\vphi)\bu+\mathbf J)-\Div(\nu(\vphi)D\bu)+\nabla \pi=\mu\nabla \vphi,\quad\text{ in }\Omega\times(0,\infty)
    \\& \label{A1AGG1}
    \Div\bu=0,\quad\text{ in }\Omega\times(0,\infty)
       \\
       \label{A1AGG2}&\partial_t\vphi+\bu\cdot\nabla\vphi-\Div(m(\vphi)\nabla \mu)=0,\quad\text{ in }\Omega\times(0,\infty)\\&
        \mu=-\Delta\vphi+f'(\vphi),\quad\text{ in }\Omega\times(0,\infty)\\&
        \label{A1AGG3}
        \partial_n\vphi=0,\quad\text{ on }\partial\Omega\times(0,\infty)\\&
        \label{A1AGG4}
        \partial_n\mu=0,\quad\text{ on }\partial\Omega\times(0,\infty)\\&
        \label{A1AGG5}
       \bu=\mathbf 0,\quad\text{ on }\partial\Omega\times(0,\infty)\\&
       \label{A1AGG6}
        \vphi(0)=\vphi_0\quad\text{ in }\Omega\\&
         \bu(0)=\bu_0\quad\text{ in }\Omega.\label{AfAGG}
    \end{align}
Here, $\bu$ is the volume averaged velocity, $\vphi$ is the difference of the volume fractions, $\rho(s)=\frac{1+s}2\rho_1+\frac{1-s}2\rho_2$ for $s\in[-1,1]$, where $\rho_j>0$, $j=1,2$, are
the densities of the two fluids. Moreover, $\mathbf{J}$ is a relative flux, related to the diffusion of the components, which is defined by
$$
\mathbf{J}:=-\frac{\rho_1-\rho_2}2m(\vphi)\nabla\mu.
$$
We point out that, the convergence to a single equilibrium of a solution to problem \eqref{A1AGG}-\eqref{AfAGG} is an open issue as well as the same result for problem \eqref{A1}-\eqref{Af} in dimension three.

The novelties of our approach are better understood if we recall the nowadays classical approach developed in \cite{AW}. The main ingredient is the regularization of a weak solution in finite time. In particular, the order parameter $\vphi$ belongs, from some positive time $\tau$ on, to $L^\infty(\tau,\infty;H^2(\Omega))$. This allows us to consider the following $\omega$-limit set
\begin{align*}
\omega(\varphi)=\{\widetilde{\varphi}\in  H^r(\Omega):\exists\, t_n\to \infty \text{ such that }\varphi(t_n)\to \widetilde{\varphi}\text{ in }H^{r}(\Omega)\},
\end{align*}
for some $r\in( \tfrac d2,2)$, which is of course nonempty by the compact embedding $H^2(\Omega)\hookrightarrow\hookrightarrow H^r(\Omega)$, and thus the trajectories are precompact in $H^r(\Omega)$. Also, this immediately gives
\begin{align}
\lim_{t\to\infty}\dist_{H^r(\Omega)}(\vphi(t),\omega(\vphi))=0.
    \label{fundam}
\end{align}
The further fundamental step is to prove that the solution is strictly separated from pure phases for $t$ large enough. Thanks to the following facts, namely, $\omega(\vphi)$ is compactly embedded in $L^\infty(\Omega)$, due to the choice of $r$, and is uniformly strictly separated from pure phases, the authors prove that there exist $T>0$ and $\delta\in(0,1)$ such that
\begin{align}
\sup_{t\geq T}\norm{\vphi(t)}_{L^\infty(\Omega)}\leq 1-\delta.
    \label{sepa1}
\end{align}
Then, after showing that any element in the $\omega$-limit set is also a stationary point for the Cahn--Hilliard equation, using \eqref{sepa1} together with \eqref{fundam}, the \L ojasiewicz-Simon inequality can be applied
and, exploiting the energy identity
\begin{align}
\label{energyid}
\frac{d}{dt}E(\vphi)+\norm{\nabla \mu}^2_{\bLd}=0,
\end{align}
conclude that $\nabla \mu\in L^1(\widetilde T,\infty;\bLd)$. Thus $\partial_t\vphi\in L^1(\widetilde T,\infty;\Hu')$, for some $\widetilde T>0$ sufficiently large. This entails that $\omega(\vphi)=\{\vphi_\infty\}$ and, by relative compactness,
\begin{align}
\norm{\vphi(t)-\vphi_\infty}_{H^{r}(\Omega)}\to0,\text{ as }t\to\infty.
    \label{conv0}
\end{align}
This argument is very robust and can be applied also to more complicated systems (see, e.g., \cite{AGGG,GPCAC,GGPS, AGGmulti,HeWu}). The main drawback is that it requires that each weak solution regularizes in finite time, which might not be the case in many situations of physical interest, as for the Cahn--Hilliard equation with non-degenerate mobility, especially when $d=3$. In these situations the convergence to a unique equilibrium has been an interesting and challenging open problem since \cite{AW}. Moreover, when the Navier--Stokes system is involved then the regularization issue is technically demanding even in the case of constant mobility (see \cite{Abels}).

The goal of this contribution is to propose a novel argument to prove the convergence to a single equilibrium of a solution to problem \eqref{A1}-\eqref{Af} which is more flexible than the one devised in \cite{AW} and widely applicable to other models, as it only requires the existence of a global weak solution satisfying an energy inequality.

Let us describe our strategy. We consider a global weak solution $\vphi$ to \eqref{A1}-\eqref{Af} such that,
for any $t\geq 0$ and almost any $0\leq s\leq t$, with $s=0$ included, there holds
    	\begin{align}
E(\vphi(t))+\int_s^t\int_\Omega m(\vphi(\tau))\normmm{\nabla\mu(\tau)}^2\dx\d\tau\leq E(\vphi(s)).\label{energyineqA}
    	\end{align}
First, we define a weaker notion of $\omega$-limit, say
\begin{align*}
\omega(\varphi)=\{\widetilde{\varphi}\in  H^1(\Omega):\exists t_n\to \infty \text{ such that }\varphi(t_n)\rightharpoonup  \widetilde{\varphi}\ \text{   weakly in }H^{1}(\Omega)\},
\end{align*}
which is nonempty, since by the energy inequality $\vphi\in L^\infty(0,\infty;H^1(\Omega))$. Then we show that actually the trajectories are precompact in $H^1(\Omega)$, a nontrivial fact on account of the low regularity of the parameter $\vphi$. Therefore, the weak convergence in the definition of $\omega(\vphi)$ can actually be replaced by the strong convergence in $\Hu$. We also prove that $\omega(\vphi)$ contains only equilibrium points for the Cahn--Hilliard equation. Then the main technical novelty is to find an alternative to the asymptotic strict separation \eqref{sepa1}, which apparently does not hold in this case, since we can only prove that
\begin{align}
\lim_{t\to\infty}\dist_{H^1(\Omega)}(\vphi(t),\omega(\vphi))=0,
    \label{fundam1}
\end{align}
and $H^1(\Omega)$ clearly does not embed into $L^\infty(\Omega)$ for $d=2,3$. The leading idea is that we can weaken the condition \eqref{sepa1} by means of a simple convergence of the Lebesgue measure of a measurable set. Namely, since we can prove that $\omega(\vphi)$ is uniformly separated from pure phases, using \eqref{fundam1} (actually the $L^2(\Omega)$-distance is enough for this purpose) we show that there exists $\delta_1>0$ such that
\begin{align}
\normmm{\{x\in\Omega:\ \normmm{\vphi(x,t)}\geq 1-\delta_1\}}\to 0,\quad\text{ as }t\to \infty.
    \label{measure}
\end{align}
Then, we consider two complementary sets of times, that is, given $T>0$ sufficiently large and $M>0$, we consider
$$
A_M(T):=\{t\in[T,\infty):\ \norm{\nabla\mu(t)}_{\bLd}\leq M\}.
$$
This is the set of ``good" times for which the dissipative term $\norm{\nabla\mu(t)}_{\bLd}$ appearing in the energy inequality \eqref{energyineqA} is uniformly controlled by the constant $M$. On the other hand, its complement $[T,\infty)\setminus A_M(T)$ is the set of ``bad" times when the dissipative term is not under control. Interestingly, a similar idea of splitting the dissipative term in ``good" and ``bad" times was also exploited in other contexts, like, for instance, to prove the weak-strong uniqueness of the Mullins-Sekerka flow, which corresponds to the sharp interface limit flow for the Cahn--Hilliard equation with non-degenerate mobility (see \cite{FHLS}).

We then focalize the attention on the ``good" times. Exploiting \eqref{measure}, we show that, given $T>0$ sufficiently large, for all $t\in A_M(T)$ a uniform strict separation property from pure phases holds. In other words,  there exists $\delta>0$ such that
\begin{align}
\sup_{t\in A_M(T)}\norm{\vphi(t)}_{L^\infty(\Omega)}\leq 1-\delta.
    \label{sepa2}
\end{align}
The proof of this result requires a De Giorgi's iterations approach (as first introduced, in the case of the Cahn-Hilliard equation, in \cite{GalP}) and crucially exploits the minimal regularity of a weak solution. This is the most delicate step of our contribution.

 Using \eqref{fundam1} and \eqref{sepa2}, we are then able to perform the \L ojasiewicz-Simon approach, this time simply making use of the energy inequality \eqref{energyineqA}, by applying a lemma firstly developed in \cite{FS} (see Lemma \ref{Feireisl} in the Appendix). This allows us to prove that $\nabla \mu\in L^1(\widetilde T,\infty;\bLd)$, for some $\widetilde T>0$ sufficiently large. The idea behind this argument is that the \L ojasiewicz-Simon inequality holds, for $T$ sufficiently large, over all the ``good" times in $A_M(T)$, whereas the set of ``bad" times can be suitably controlled.
Finally, we conclude the argument, proving that $\omega(\vphi)=\{\vphi_\infty\}$, and, by relative compactness,
\begin{align}
\norm{\vphi(t)-\vphi_\infty}_{H^{1}(\Omega)}\to0,\text{ as }t\to\infty.
    \label{conv2}
\end{align}


We expect that this approach can also be applied to similar equations, like, for instance, the Allen-Cahn equation with non-constant mobility, as it will be shown in a forthcoming contribution.
\medskip

The plan of the paper is the following. The next section is devoted to some preliminaries which include the statements of the existence of a global weak solution to problem \eqref{A1}-\eqref{Af} and problem \eqref{A1AGG}-\eqref{AfAGG}. Then, in Section \ref{main}, we state our main results. Two fundamental lemmas are proven in Section \ref{secconvaaa} and Section \ref{sec_prooftwoparts}, respectively. The proof
of the main result about the Cahn--Hilliard equation is given in Section \ref{sec:proofLoja}, while the one related to the AGG model is contained in Section \ref{sec:proofLoja}. Two known technical lemmas are reported in Appendix for the reader's convenience.
    \section{Preliminaries}
    \label{pre}
    \subsection{Notation and functional setting}
Here we introduce some notation along with the functional spaces which will be used in the sequel.

 \begin{enumerate}[label=\textnormal{(N\arabic*)},leftmargin=*]
 	
 	\item \textbf{Notation for general Banach spaces.}
 	For any (real) normed space $X$
    , we denote its norm by $\|\cdot\|_X$,
 	its {dual space by $X'$}.
 	If $X$ is a Hilbert space, we write $(\cdot,\cdot)_X$ to denote the corresponding inner product.
 	Moreover, the corresponding spaces of vector-valued or matrix-valued functions with each component in $X$ are denoted by $\mathbf{X}$.
 	
 	\item \textbf{Lebesgue and Sobolev spaces.}
 	Assume $\Omega$ to be a bounded domain in $\R^d$, $d=2,3$ of class $C^2$.
 	For $1 \leq p \leq \infty$ and $k \in \N$, the classical Lebesgue and Sobolev spaces defined on $\Omega$ are denoted by $L^p(\Omega)$ and $W^{k,p}(\Omega)$, and their standard norms are denoted by $\|\cdot\|_{L^p(\Omega)}$ and $\|\cdot\|_{W^{k,p}(\Omega)}$, respectively.
 	In the case $p = 2$, we set $H^k(\Omega) = W^{k,2}(\Omega)$. Note that the $L^2(\Omega)$ inner product is simply denoted by $(\cdot,\cdot)$.
 	Also, for any interval $I\subset\R$, any Banach space $X$, $1 \leq p \leq \infty$ and $k \in \N$, we write $L^p(I;X)$, $W^{k,p}(I;X)$ and $H^{k}(I;X) = W^{k,2}(I;X)$ to denote the Lebesgue and Sobolev spaces of functions with values in $X$. The canonical norms are indicated by $\|\cdot\|_{L^p(I;X)}$, $\|\cdot\|_{W^{k,p}(I;X)}$ and $\|\cdot\|_{H^k(I;X)}$, respectively.
 	    We additionally define
\begin{align*}
    L^p_\mathrm{loc}(I;X)
    &:=
    \big\{
        u:I\to X \,\big\vert\, u \in L^p(J;X) \;\text{for every compact interval $J\subset I$}
    \big\}
    \\[1ex]
    L^p_\mathrm{uloc}(I;X)
    &:=
    \left\{ u:I\to X \,\middle|\,
    \begin{aligned}
    &u \in L^p_\mathrm{loc}(I;X) \;\text{and}\; \exists\, C>0\; \sup_{t\in I}\|u\|_{L^p(t,t+1;X)} \le C
    \end{aligned}
    \right\}.
\end{align*}
The spaces $W^{k,p}_\mathrm{loc}(I;X)$, $H^k_\mathrm{loc}(I;X)$, $W^{k,p}_\mathrm{uloc}(I;X)$, $H^k_\mathrm{uloc}(I;X)$ are defined analogously.
 	\item \textbf{Spaces of continuous functions.}
 	For any interval $I\subset\R$ and any Banach space $X$, $C(I;X)$ denotes the space of continuous functions mapping from $I$ to $X$ and $BC(I;X)$ denotes the space of bounded functions in $C(I;X)$. Moreover, $C_\mathrm{w}(I;X)$ denotes the space of functions mapping from $I$ to $X$, which are continuous on $I$ with respect to the weak topology of $X$, and $BC_\mathrm{w}(I;X)$ denotes the space of bounded functions in $C_\mathrm{w}(I;X)$. Then, we denote by $C^\gamma(I;X)$, $\gamma\in(0,1]$, the space of $\gamma$-H\"{o}lder (Lipschitz, if $\gamma=1$) continuous functions with values in $X$. In addition, $C_0^k(I;X)$ stands for the space of $k$-continuously differentiable functions with compact support mapping $I$ into $X$.
 	\item \textbf{Spaces of functions with zero integral mean.}
If $v\in H^1(\Omega)'$, then its generalized spatial mean is defined as
\begin{equation*}
    \overline{v}:= \frac{\langle v,1 \rangle_{\Hu',\Hu}}{\normmm{\Omega}},
\end{equation*}%
where $|\Omega |$ stands for the $d$-dimensional Lebesgue measure of
$\Omega$. Clearly, if $v\in L^1(\Omega)$, then this spatial mean becomes the classical integral average
\begin{align*}
    \overline v=\frac{\int_\Omega v\dx}{\normmm{\Omega}}.
\end{align*}
Using this definition, we introduce the following Hilbert spaces:
\begin{align*}
    L^2_{(0)}(\Omega) &:= \big\{ u\in L^2(\Omega) \,:\, \overline u = 0 \big\} \subset L^2(\Omega),\\
    H^1_{(0)}(\Omega) &:= \big\{ u\in H^1(\Omega) \,:\, \overline u = 0 \big\} \subset H^1(\Omega),\\
     H^{1}_{(0)}(\Omega)' &:= \big\{ u\in H^1(\Omega)' \,:\, \overline u = 0 \big\} \subset H^1(\Omega)'.
\end{align*}
Additionally, for $k\in \R$, we define the affine space
$$
H^1_{(k)}(\Omega)=H^1_{(0)}(\Omega)+k.
$$

 	\item \textbf{Spaces of divergence-free functions.}
 We define the closed linear subspaces
 	\begin{align*}
 		\mathbf L^p_\sigma(\Omega)
 		&:=\overline{\{\mathbf{u}\in \mathbf{C}^\infty_0(\Omega) \,\big\vert\, \operatorname{div}\ \mathbf{u}=0\}}^{\mathbf{L}^p(\Omega)}
 		\subset \mathbf L^p(\Omega), \quad p\in[2,\infty),\\
 		\mathbf H^1_\sigma(\Omega)&:= \mathbf L^2_\sigma(\Omega) \cap \mathbf H^1(\Omega).
 	\end{align*}
In both cases, Korn's inequality yields
\begin{equation}
\Vert \mathbf{u}\Vert \leq \sqrt{2}\Vert D\mathbf{u}%
\Vert\leq \sqrt{2}\Vert \nabla \mathbf{u}\Vert
\quad \text{ for all } \mathbf{u}\in \mathbf H^1_\sigma (\Omega).
\label{korn}
\end{equation}
As a trivial consequence, $\|\nabla\cdot\|$ is a norm on $\mathbf H^1_\sigma(\Omega)$ that is equivalent to the standard norm $\|\cdot\|_{\mathbf{H}^1(\Omega)}$.
 \end{enumerate}

    \subsection{Main assumptions}
    We enumerate here all the assumptions that are needed to establish the results of this contribution. We start with the ones concerning the Cahn--Hilliard equation with non-degenerate mobility:
\begin{enumerate}[label=\textnormal{(A\arabic*)},leftmargin=*]
    \item\label{ASS:0} The non-degenerate mobility function $m(\cdot)$ satisfies the following conditions
\begin{align*}
    m\in C([-1,1]),\quad 0<m_*\leq m(s)\leq m^*,\quad \forall s\in[-1,1].
\end{align*}
    \item \label{ASS:S1}
    The potential $f:[-1,1]\to \R$ can be written as follows
    \begin{equation*}
        f(s)=F(s)-\frac{\theta_0}{2}s^2 \quad\text{for all $s\in [-1,1]$}
    \end{equation*}
    with a given constant $\theta_0>0$, where $F\in C([-1,1])\cap C^{2}(-1,1)$ is such that
    \begin{equation*}
    \lim_{r\rightarrow -1}F^{\prime }(r)=-\infty ,
    \quad \lim_{r\rightarrow 1}F^{\prime }(r)=+\infty ,
    \quad F^{\prime \prime }(s)\geq {\theta},
    \quad F'(0)=0
    \end{equation*}
    for all $s\in (-1,1)$ and a prescribed constant $\theta\in(0,\theta_0)$.
    Without loss of generality, we further assume $F(0)=0$ and $F'(0)=0$.
    In particular, this means that $F(s)\geq 0$ for all $s\in [-1,1]$.
    For the sake of convenience, we extend $f$ and $F$ onto $\R\setminus[-1,1]$ by defining
    $f(s):=+\infty $ and $F(s):=+\infty $ for all $s\in\R\setminus [-1,1]$.

\end{enumerate}
\begin{remark}\label{REM:LOG}
Note that \eqref{FHP} can be written as 
\begin{equation}
    f(s)=F_\mathrm{log}(s)-\frac{\theta_0}{2}s^2 \quad\text{for all $s\in [-1,1]$},
    \label{f:LOG}
\end{equation}
with $F_\mathrm{log}(\pm 1) = \theta\ln(2)$ and
\begin{equation}
    F_\mathrm{log}(s)=\frac{\theta}{2}((1+s)\text{ln}(1+s)+(1-s)\text{ln}(1-s))
    \quad\text{for all $s\in(-1,1)$}.
\label{F:LOG}
\end{equation}
Thus, it satisfies all assumptions \ref{ASS:S1}.
\end{remark}
The results on problem \eqref{A1AGG}-\eqref{AfAGG} need the further assumptions:
\begin{enumerate}[label=\textnormal{(H\arabic*)},leftmargin=*]
        \item \label{ASS:Viscosity} The viscosity $\nu\in W^{1,\infty}(\R)$ satisfies
    \begin{align*}
        0 < \nu_* \leq \nu(s) \leq \nu^*\qquad \text{ for all }s\in\R,
    \end{align*}
    for some positive constants $\nu_*,\nu^*\in \R$.
    \item \label{ASS:3} The density $\rho$ is such that
    \begin{align*}
        \rho(s):=\frac{1+s}2\rho_1+\frac{1-s}2\rho_2,
    \end{align*}
    for $\rho_1,\rho_2>0$, and we set
    \begin{align*}
        \rho_*:=\min\{\rho_1,\rho_2\}>0,\quad \rho^*:=\max\{\rho_1,\rho_2\}>0.
    \end{align*}
\end{enumerate}
    \subsection{Existence of a global weak solution}
Here we state a result concerning the existence of weak solutions to problem \eqref{A1}-\eqref{Af}.
This notion of solution is the one we need to prove the convergence to a unique equilibrium.
The main result on the existence of such solutions dates back to Barrett and Blowey \cite{BB}. In dimension two, some stronger results are known (see \cite{BB,CGGG}), provided that $m$ is
smoother, but our aim here is to highlight the \textit{minimal }requirements to study the longtime behavior of trajectories and this is the only available result so far.
    \begin{theorem}
[Existence of global weak solutions \cite{BB}]
 \label{weak}
  Let $\Omega\subset R^d$, $d=2,3$, be a bounded domain of class $C^2$, and let assumptions \ref{ASS:0}-\ref{ASS:S1} be satisfied. If $\vphi_0 \in H^1(\Omega)$ is such that $\normmm{\vphi_0}\leq 1$ and $\overline\vphi_0\in(-1,1)$, then
there exists a global weak solution $(\vphi,\mu)$ to problem \eqref{A1}-\eqref{Af}. This means that
  \begin{align*}
     & \vphi\in BC([0,\infty);H^1(\Omega))\cap L^2_{uloc}([0,\infty);H^2(\Omega)),\\&
     \vphi\in L^\infty(\Omega\times(0,\infty)):\ \normmm{\vphi(x,t)}<1,\quad\text{ for a.a. }(x,t)\in \Omega\times(0,\infty),\\&
     \partial_t\vphi\in L^2(0,\infty;H^{1}_{(0)}(\Omega)'),\\&
     \mu\in L^2_{uloc}([0,\infty);H^1(\Omega)),
  \end{align*}
    and
    \begin{align*}
        &\langle \partial_t\vphi,v\rangle_{\Hu',\Hu}+(m(\vphi)\nabla \mu,\nabla v)=0,\quad \forall v\in H^1(\Omega), \text{ for a.a. }t\geq 0,\\&
        \mu=-\Delta\vphi+f'(\vphi)\quad\text{ a.e. in }\Omega\times(0,\infty),
    \end{align*}
    together with $\partial_{\bn}\vphi=0$ on $\partial\Omega\times(0,\infty)$.
    Additionally, for any $t\geq 0$ and almost any $s \in [0,t]$, with $s=0$ included, it holds
    	\begin{align}
    		E(\vphi(t))+\int_s^t\int_\Omega m(\vphi(\tau))\normmm{\nabla\mu(\tau)}^2\dx\d\tau\leq E(\vphi(s)),\label{energyineq}
    	\end{align}
        where
       \begin{align*}
           E(v)=\frac12\int_\Omega \normmm{\nabla v}^2\dx+\int_\Omega f(v)\dx,
       \end{align*}
       for any $v\in H^1(\Omega)$ such that $\normmm{v}\leq 1$ in $\Omega$.
    \end{theorem}
    \begin{remark}
     To be precise the energy inequality \eqref{energyineq} can be proven to be an identity (see, for instance, \cite[Lemma 2.4]{Schimperna}). Nevertheless, in this contribution we aim at stressing the fact that our proof requires nothing more than an energy inequality.  \label{schimp}
    \end{remark}
We now state here an existence result concerning problem \eqref{A1AGG}-\eqref{AfAGG} (see \cite{ADG}). We have
    \begin{theorem}
[Weak existence of global solutions to AGG system \cite{ADG}]
 \label{weakAGG}
  Let $\Omega\subset R^d$, $d=2,3$, be a bounded domain of class $C^2$, and let assumptions \ref{ASS:0}-\ref{ASS:S1} and \ref{ASS:Viscosity}-\ref{ASS:3} hold.
  If $\bv_0\in \bLds$ and $\vphi_0 \in H^1(\Omega)$ is such that $\normmm{\vphi_0}\leq 1$ and $\overline\vphi_0\in(-1,1)$, then
there exists a global weak solution $(\bu,\vphi,\mu)$ to problem \eqref{A1AGG}-\eqref{AfAGG}. This means that
  \begin{align*}
  &\bu\in BC_w([0,\infty);\bLds)\cap L^2(0,\infty;\mathbf H^1_\sigma(\Omega)),\\&
      \vphi\in BC_w([0,\infty);H^1(\Omega))\cap L^2_{uloc}([0,\infty);H^2(\Omega)),\\&
     \vphi\in L^\infty(\Omega\times(0,\infty)):\ \normmm{\vphi(x,t)}<1,\quad\text{ for a.a. }(x,t)\in \Omega\times(0,\infty),\\&
     \partial_t\vphi\in L^2(0,\infty;H^{1}_{(0)}(\Omega)'),\\&
     \mu\in L^2_{uloc}([0,\infty);H^1(\Omega)),
  \end{align*}
    and, for any $T>0$,
    \begin{align}
& \int_0^T\left(-(\rho (\vphi )\mathbf{u}),\partial_t\mathbf{w})-\left( \rho (\vphi )\mathbf{u}\otimes \mathbf{u}%
,D\mathbf{w}\right) +\left( \nu (\vphi )D\mathbf{u},D\mathbf{w}\right)
-\left( \mathbf{u},\left( \mathbf{J}\cdot \nabla \right) \mathbf{w}\right)\right)\d\tau
\nonumber\\&=-\int_0^T\left( \vphi \nabla \mu ,\mathbf{w}\right)\d\tau , \text{ for any }\mathbf{w}\in [C^\infty_0(\Omega\times(0,T)]^d, \text{ with }\Div\bw=0,\\
& \langle \partial _{t}\vphi ,v\rangle _{\Hu',\Hu}-\left( \vphi \,%
\mathbf{u},\nabla v\right) +\left( m(\vphi )\nabla \mu ,\nabla v\right) =0,\quad  \forall v\in H^1(\Omega)\text{ and a.a. }t\in(0,T),  \label{CHeq}
\\&
        \mu=-\Delta\vphi+f'(\vphi)\quad\text{ a.e. in }\Omega\times(0,T),
\label{weak-NS}
\end{align}%
together with $\partial_{\bn}\vphi=0$ on $\partial\Omega\times(0,\infty)$.
    Additionally, the following energy inequalities hold: for any $t\geq 0$ and almost any $s \in [0,t]$, with $s=0$ included, it holds
    	    	\begin{align}
&\nonumber\frac12\int_\Omega \rho(\vphi(t))\normmm{\bu(t)}^2\dx+\int_s^t\int_\Omega\nu(\vphi(\tau)\normmm{D\bu(\tau)}^2\dx\d\tau\\&\leq \frac12\int_\Omega \rho(\vphi(s))\normmm{\bu(s)}^2\dx-\int_s^t\int_\Omega \bu(\tau)\cdot\nabla\mu(\tau)\vphi(\tau)\dx\d\tau,\label{energyineq3}
    	\end{align}
    as well as
\begin{align}
    		&E(\vphi(t))+\int_s^t\int_\Omega m(\vphi(\tau))\normmm{\nabla\mu(\tau)}^2\dx\d\tau\leq E(\vphi(s))+\int_s^t\int_\Omega \bu(\tau)\cdot\nabla\mu(\tau)\vphi(\tau)\dx\d\tau.\label{energyineq3b}
    	\end{align}
    \end{theorem}
\begin{remark}
To be precise, the existence result proven in \cite{ADG} only shows the validity of the total energy inequality, that is, for any $t\geq 0$ and almost any $s\in [0,t]$, with $s=0$ included, it holds
\begin{align}
    		&\nonumber E_{tot}(\vphi(t),\bu(t))+\int_s^t\int_\Omega\nu(\vphi(\tau)\normmm{D\bu(\tau)}^2\dx\d\tau\\&+\int_s^t\int_\Omega m(\vphi(\tau))\normmm{\nabla\mu(\tau)}^2\dx\d\tau\leq E_{tot}(\vphi(s),\bu(s)),\label{energyineq2}
    	\end{align}
        where
        \begin{align}
        E_{tot}(\bu,\vphi):=\frac12\int_\Omega \rho(\vphi)\normmm{\bu}^2\dx+E(\vphi).
            \label{Etot}
        \end{align}
Note that this estimate can be obtained by summing the two energy inequalities \eqref{energyineq3} and \eqref{energyineq3b}. Nevertheless, by studying the approximation scheme adopted in \cite{ADG} to construct the weak solution (see in particular \cite[Eqs. (4.8)-(4.9)]{ADG}), and then passing to the limit in the approximating parameter $k\to \infty$, we can easily see that \eqref{energyineq3} and \eqref{energyineq3b} actually hold separately. Indeed, it is enough to argue as in \cite[Section 5.3]{ADG}, but considering the two inequalities in a separate fashion, without adding them up. A similar energy splitting is used, for instance, in the definition of weak solution \cite[Definition 3.4]{AGP1}, and we also refer to \cite[Proof of Theorem 3.8, 5.7.7]{AGP1} for a proof on how to pass to the limit in the two inequalities.

\end{remark}

\begin{remark}
        In both the problems the conservation of mass holds, that is, $\overline{\vphi}(t)=\overline{\varphi_0}$ for all $t\geq 0$.
\end{remark}

\begin{remark}
        In \cite{ADG} the authors actually prove the result for the case of mobility $m\in C^1([-1,1])$, but the case $m\in C([-1,1])$ can be easily obtained by an approximation argument.
\end{remark}

\section{Main results}
\label{main}
  \subsection{The Cahn--Hilliard equation with non-degenerate mobility}\label{mainCH}
	We discuss the longtime behavior of each (weak) trajectory. Let us consider the set of admissible initial data:
\begin{align}
\mathcal{H}_k:=\left\{\varphi\in H^1(\Omega): \Vert\varphi\Vert_{L^\infty(\Omega)}\leq 1,\quad \vert\overline{\varphi}\vert= k \right\}, \label{Hk}
\end{align}
with $k\in[0,1)$, and fix an initial datum $\varphi_0\in \mathcal{H}_k$. Let then $\varphi$ be a weak global-in-time solution departing from $\varphi_0$, which might be nonunique, whose existence is ensured by Lemma \ref{weak}. We introduce the (weak) $\omega$-limit set associated to $\varphi$, i.e.,
\begin{align*}
\omega(\varphi)=\{\widetilde{\varphi}\in  \mathcal{H}_k:\exists t_n\to \infty \text{ such that }\varphi(t_n)\rightharpoonup \widetilde{\varphi}\ \text{ weakly in }H^1(\Omega)\}.
\end{align*}
Note that $\varphi$ is uniformly bounded in $H^1(\Omega)$, which is a reflexive space. Therefore $\omega(\vphi)$ is non-empty. We now characterize the set $\omega(\varphi)$, showing that it is composed by equilibrium points
according to the following definition.
\begin{definition}
	$\varphi_\infty$ is an equilibrium point to the Cahn--Hilliard equation with non-degenerate mobility \eqref{A1}-\eqref{Af} if $\varphi_\infty\in \mathcal{H}_k\cap H^2(\Omega)$ satisfies the stationary Cahn--Hilliard equation
	\begin{align}
		-\Delta \varphi_\infty+f^\prime(\varphi_\infty)=\mu_\infty,\quad \text{ a.e. in }\Omega,
		\label{conv1t}
	\end{align}
	together with $\partial_{\bn}\varphi_\infty=0$ on $\partial\Omega$, where $\mu_\infty\in \R$ is a real constant.
\end{definition}
\begin{remark}
    Note that, given $\mu_\infty\in \R$, solutions to \eqref{conv1t} do exist (see, e.g., \cite{AW, GGPS}), but they might form a set with the power of continuum.
\end{remark}
If we introduce the set of all the stationary points of the Cahn--Hilliard equation:
$$
\mathcal{S}:=\left\{\varphi_\infty\in \mathcal{H}_k\cap H^2(\Omega): \varphi_\infty\text{ satisfies }\eqref{conv1t}\right\},
$$
we can easily prove that $\omega(\vphi)\subset \mathcal{S}$. In particular, in Section \ref{secconvaaa} we will prove the following
\begin{lemma}
				\label{convaaa}
				Let the assumptions of Theorem \ref{weak} hold. We have
				$$
				    \omega(\vphi)\subset \mathcal{S}.
				$$
				Moreover, $\omega(\vphi)$ is bounded in $ H^2(\Omega)$, and there exists $\delta_1>0$ such that
				\begin{align}
				\| \vphi_\infty\|_{L^\infty(\Omega)}\leq 1-2\delta_1,\quad \forall \: \vphi_\infty\in \omega(\vphi).
				\label{sepaglobal}
                \end{align}
In conclusion, the trajectories of $\vphi(\cdot)$ are precompact in $H^1(\Omega)$, so that $\omega(\vphi)$ is compact in $H^1(\Omega)$, it holds the characterization
                \begin{align}
\omega(\varphi)=\{\widetilde{\varphi}\in  \mathcal{H}_k:\exists t_n\to \infty \text{ such that }\varphi(t_n)\to \widetilde{\varphi}\text{ in }H^1(\Omega)\},\label{omegal}
\end{align}
 and it holds
\begin{align}
	\lim_{t\to\infty}\dist_{ H^1(\Omega)}(\varphi(t),\omega(\vphi))=0.\label{convergence}
\end{align}
			\end{lemma}
  Using this lemma, we can prove the following properties, which are crucial to show that $\omega(\vphi)$ is indeed a singleton (see Section \ref{sec_prooftwoparts} for the proof).
\begin{lemma}\label{twoparts}
    Let the assumptions of Theorem \ref{weak} hold. Given the set
    \begin{align}
    	A_\delta(t) &:= \{x\in\Omega:\; |\varphi(x,t)|\geq 1-{\delta_1}\},\quad t\geq0 ,\label{Adelta}
    \end{align}
    it holds
    \begin{align}
    	\lim_{t\to\infty }\normmm{A_\delta(t)}\to 0,
    	\label{conerg}
    \end{align}
    where $\delta_1>0$ is given in \eqref{sepaglobal}.
   Then, for any $M>0$ there exists $\delta\in(0,\delta_1)$ and $T_S>0$ such that
    \begin{align}
        \sup_{t\in A_M(T_S)}\norm{\varphi(t)}_{L^\infty(\Omega)}\leq 1-\delta, \label{asympt}
    \end{align}
    where $$
    A_M(T_S):=\{t\geq T_S:\ \norm{\nabla\mu(t)}_{\bLd}\leq  M\}.
    $$
\end{lemma}
As a consequence of this fundamental lemma, in Section \ref{sec:proofLoja} we can finally prove that the $\omega$-limit is formed by a unique element, if the potential is analytic in $(-1,1)$, namely,
\begin{theorem}\label{uniqueeq}
	Let the assumptions of Theorem \ref{weak} hold and suppose additionally that $F$ is real analytic in $(-1,1)$. Then any global weak solution $\vphi$ given by Theorem \ref{weak}, departing from the initial datum $\vphi_0\in
    \mathcal{H}_k$, converges to a single equilibrium point $\vphi_\infty\in \mathcal S$, i.e., $\omega(\vphi)=\{\vphi_\infty\}$. In particular, it holds
	\begin{align} \lim_{t\to \infty}\Vert\vphi(t)-\vphi_\infty\Vert_{H^1(\Omega)}=0.
		\label{equil}
	\end{align}

 \end{theorem}
 \begin{remark}
     In dimension two, if $m\in C^2([-1,1])$, the uniqueness of weak solutions hold (see \cite{CGGG}). Thus, our result gives the convergence of each (unique) trajectory to a unique equilibrium point $\vphi_\infty\in \mathcal S$. However, compared to the asymptotic convergence result given in \cite{CGGG}, here we can weaken the assumptions on $m$ (see \ref{ASS:0}) as well as on $F$, which is here only required to be singular at pure phases (see \ref{ASS:S1} and cf. \cite[Rem.1.3]{CGGG}).
 \end{remark}
 \subsection{The Abels-Garcke-Gr\"un system with non-degenerate mobility}\label{secAGG}
 Here we show the robustness of the method by proving the convergence to a single equilibrium of the weak solutions to problem \eqref{A1AGG}-\eqref{AfAGG}. This is an open issue, due to the lack of instantaneous (and even asymptotic) regularization, especially in three dimensions. However, our argument works as it only exploits the existence of a global weak solution satisfying an energy inequality.

 Let us fix the initial data $(\bu_0,\varphi_0)\in \mathcal{H}_k\times \bLds$, where $\mathcal H_k$ is introduced in \eqref{Hk}. Consider then a global weak solution$(\bu,\varphi)$ originating from $(\bu_0,\varphi_0)$ whose existence is ensured by Theorem \ref{weakAGG}.

 First, only relying on the energy inequalities \eqref{energyineq2} and \eqref{energyineq3}, we can prove the following lemma (see Section \ref{secproofconvu} for its proof):
 \begin{lemma}\label{lemmauconv}
     Let the assumptions of Theorem \ref{weakAGG} hold. Then, for any fixed $(\bu_0,\vphi_0)\in \mathcal H_k\times \bLds$, a corresponding global weak solution $(\bu,\vphi)$ given by Theorem \ref{weakAGG} is such that
     \begin{align}
\lim_{t\to\infty}\norm{\bu(t)}_{\bLds}=0,
         \label{uconv}
     \end{align}
     i.e., $\bu(t)\to\mathbf 0$ in $\bLds$ as $t\to\infty$.
 \end{lemma}
 We then introduce the $\omega$-limit set associated to $(\bu,\varphi)$, i.e.,
\begin{equation*}
\omega(\bu,\varphi):=\{(\mathbf 0,\widetilde{\varphi})\in  \{\mathbf 0\}\times \mathcal{H}_k:\exists\, t_n\to \infty \text{ s.t. }\varphi(t_n)\rightharpoonup \widetilde{\varphi}\text{ in }H^1(\Omega),\;\bu(t_n)\to\mathbf 0\text{ in }\bLds\}.
\end{equation*}
Observe that $\varphi$ is uniformly bounded in $H^1(\Omega)$ and thus, since by \eqref{uconv} we have $\bu(t)\to 0$ in $\bLds$ as $t\to\infty$, then $\omega(\bu,\vphi)$ is non-empty. We can show now that $\omega(\bu,\varphi)$ contains only by equilibrium points which are defined as follows (cf. \eqref{conv1t}).
\begin{definition}
	$(\mathbf 0,\varphi_\infty)$ is an equilibrium point to problem \eqref{A1AGG}-\eqref{AfAGG} if $\varphi_\infty\in \mathcal{H}_k\cap H^2(\Omega)$ satisfies the stationary Cahn--Hilliard equation
	\begin{align}
		-\Delta \varphi_\infty+f^\prime(\varphi_\infty)=\mu_\infty,\quad \text{in }\Omega,
		\label{conv1t1}
	\end{align}
	together with $\partial_{\bn}\varphi_\infty=0$ on $\partial\Omega$, where $\mu_\infty\in \R$ is a real constant.
\end{definition}
Let us now set
$$
\mathcal{S}_1:=\left\{(\mathbf0,\varphi_\infty)\in \{\mathbf 0\}\times\mathcal{H}_k\cap H^2(\Omega): \varphi_\infty\text{ satisfies }\eqref{conv1t1}\right\}.
$$
Then, we can easily prove that $\omega(\bu,\vphi)\subset \mathcal{S}_1$. Also, arguing as in Section \ref{mainCH}, we can finally prove (see Section \ref{secprooffinal}) our main result, i.e., the property that $\omega(\bu,\vphi)$ is a singleton, as long as we assume $F$ to be real analytic in $(-1,1)$. Namely we have
\begin{theorem}\label{uniqueeq1}
	Let the assumptions of Theorem \ref{weakAGG} hold. We have
				$$
				    \omega(\bu,\vphi)\subset \mathcal{S}_1.
				$$
    Furthermore $\omega(\bu,\vphi)$ is compact in $H^1(\Omega)$, bounded in $\{ \mathbf 0\} \times H^2(\Omega)$, and there exists $\delta_1>0$ such that
				\begin{align}
				\| \vphi_\infty\|_{L^\infty(\Omega)}\leq 1-2\delta_1,\quad \forall \: (\mathbf0,\vphi_\infty)\in \omega(\bu,\vphi).
				\label{sepaglobal1}
                \end{align}
                Additionally, there holds
                \begin{align}
	\lim_{t\to\infty}\dist_{ \bLds\times H^1(\Omega)}((\bu(t),\varphi(t)),\omega(\bu,\vphi))=0,\label{convergence1}
\end{align}
and, in particular,
\begin{align}\label{convu1}
                    \bu(t)\to \mathbf 0\text{ in }\bLds,\quad \text{  as }t\to\infty.
                \end{align}
    If, in addition, we suppose that $F$ is real analytic in $(-1,1)$, then any global weak solution $(\bu,\vphi)$, originating from $(\bu_0,\vphi_0)\in \bLds\times \mathcal{H}_k$, converges to a single equilibrium point $(\mathbf 0,\vphi_\infty)\in \mathcal S_1$ and $\omega(\bu,\vphi)=\{(\mathbf 0,\vphi_\infty)\}$. In particular, we have
	\begin{align} \lim_{t\to \infty}\Vert\vphi(t)-\vphi_\infty\Vert_{H^1(\Omega)}=0.
		\label{equil1}
	\end{align}
 \end{theorem}

 \section{Proof of Lemma \ref{convaaa}}\label{secconvaaa}
Let us consider a sequence $t_n\to \infty$ such that $\varphi(t_n)\rightharpoonup  \widetilde{\varphi}$ weakly in $H^1(\Omega)$, with $\widetilde\vphi\in \omega(\varphi)$. Of course by compactness, we can focalize our attention on a nonrelabeled subsequence such that $\vphi(t_n)\to\widetilde\vphi$ strongly in $\Ld$. We then define the sequence of trajectories $\varphi_n(t):=\varphi(t+t_n)$ and $\mu_n(t):=\mu(t+t_n)$. They solve
\begin{align}
      \label{eqs1}  &\langle \partial_t\vphi_n,v\rangle_{\Hu',\Hu}+(m(\vphi_n)\nabla \mu_n,\nabla v)=0,\quad \forall v\in H^1(\Omega), \text{ for a.a. }t\geq 0,\\&
        \mu_n=-\Delta\vphi_n+f'(\vphi_n)\quad\text{ a.e. in }\Omega\times(0,\infty),\label{eqs2}
    \end{align}
   together with $\partial_{\bn}\varphi_n=0$ almost everywhere on $\partial\Omega \times [0,+\infty)$,  
Thanks to Theorem \ref{weak}, recalling the energy inequality \eqref{energyineq}, we get that ${{E}}(\vphi(t_n))\leq {{E}}(\vphi_0)$ for any $n$, and thus, for any $T>0$, there exists $C(T)>0$ independent of $n$ such that
				\begin{align}
				\| \vphi_n\|_{L^\infty(0,T;H^1(\Omega))}+\| \mu_n\|_{L^2(0,T;H^1(\Omega))}\leq C(T).\label{energy}
				\end{align}
                Here the second term  can be controlled, for instance, via energy inequality and Poincaré's inequality, using the mass conservation and the well-known control (see, e.g., \cite{Kenmochi, FrigeriGrasselli})
                $$
                \int_\Omega\normmm{F'(\vphi_n)}\dx\leq C(\overline\vphi_0)\left(1+\int_\Omega F'(\vphi_n)(\vphi_n-\overline\vphi_n)\dx\right),
                $$
                which gives
                \begin{align}
                   \int_\Omega\normmm{\mu_n}\dx=\int_\Omega\normmm{f'(\vphi_n)}\dx\leq C(1+\norm{\nabla\mu}_{\bLd}).\label{ctr}
                \end{align}
Then, by comparison, we also get
\begin{align}
    \| \partial_t\vphi_n\|_{L^2(0,T;H^1(\Omega)')}\leq C(T).\label{cmp}
\end{align}
Moreover, testing (formally, but it can easily be made rigorous by suitable approximations) the equation for $\mu_n$ with $F'(\vphi_n)$ and recalling \eqref{ctr}, we get
$$
\norm{F'(\vphi_n)}_{L^2(0,T;L^2(\Omega))}\leq C(1+\norm{\nabla\mu}_{L^2(0,T;\bLd)})\leq C(T).
$$
Thus, by elliptic regularity, we find
\begin{align}
\| \vphi_n\|_{L^2(0,T;H^2(\Omega))}\leq C(T).
    \label{elliptic}
\end{align}
From estimates \eqref{energy}-\eqref{elliptic}, we deduce that there exists $\vphi^*$ (and $\mu^*$) such that, for any fixed $T>0$,
				\begin{align}
	&
        \label{est}\vphi_n\rightharpoonup \vphi^*\quad\text{weakly in } L^2(0,T;\Hu)\cap H^1(0,T;\Hu'),\\&
					\vphi_n\overset{\ast}{\rightharpoonup} \vphi^*\quad\text{weakly-$*$ in } L^\infty(\Omega\times(0,T)),
					\\&
					\vphi_n\to \vphi^*\quad \text{in } L^2(0,T;H^s(\Omega)),\, \forall\: s\in[0,2)\text{ and a.e. in }\Omega\times(0,T),\label{essenziale}\\&
					\mu_n\rightharpoonup  \mu^*\quad\text{weakly in } L^2(0,T;\Hu).
				\end{align}
				Additionally, by the Aubin-Lions Lemma we also infer that, for any $T>0$,
				\begin{align}
					\vphi_n\to \vphi^*\quad\text{strongly in }C([0,T];L^2(\Omega)).
					\label{AL}
				\end{align}
				These convergences are enough to pass to the limit in the equations \eqref{eqs1}-\eqref{eqs2}, meaning that {the limit pair } $(\vphi^*,\mu^*)$ satisfies, for any $T>0$,
				\begin{align*}
					\langle \partial_t\vphi^*,v\rangle+(m(\vphi^*)\nabla \mu^*, \nabla v)=0,\quad \forall \: v\in \Hu,\quad \text{a.e. in }(0,T),\\
					\mu^*=-\Delta\vphi^*+f'(\vphi^*),\quad\text{a.e. in }\Omega\times(0,T),
				\end{align*}
				with initial datum $\vphi^*(0)={\vphi_\infty}$ and boundary condition $\partial_{\bn}\varphi^*=0$ on $\partial\Omega$. This follows immediately from the fact that $\vphi_n(0)=\vphi(t_n){\to} {\vphi_\infty}$ strongly in $\Ld$. Furthermore, we have \[
				\lim_{n\to \infty}{{E}}(\vphi_n(t))={{E}}(\vphi^*(t))
				\] for {almost any} $t\geq 0$. By the energy inequality \eqref{energyineq}, we infer that the energy ${{E}}(\vphi(\cdot))$ is nonincreasing, thus there exists ${{E}}_\infty$ such that
                \begin{align}
				\lim_{t\to+ \infty}{{E}}(\vphi(t))={{E}}_\infty.
				\label{ene}
                \end{align}
				 Hence, {for almost any $t\geq0$}, we have
				$${{E}}(\vphi^*(t))=\lim_{n\to \infty}{E}(\vphi_n(t))=\lim_{n\to \infty}{{E}}(\vphi(t+t_n))={{E}}_\infty,$$so that $E(\vphi^*(\cdot))$ is constant in time and equal to $ E_\infty$. Passing then to the limit in the energy inequality, which is valid for each $\vphi_n$ thanks again to \eqref{energyineq}, we obtain
				
			\begin{align}\label{zeros}
				{{E}}_\infty+\int_s^{t}\int_\Omega m(\vphi^*(\tau))\normmm{\nabla\mu^*(\tau)}^2\dx \: \d \tau\leq {{E}}_\infty\quad \text{ for almost any } 0\leq s\leq t<\infty,\end{align}
			with $s=0$ included.

    Then, since by assumption $m(\cdot)\geq m_*>0$, \eqref{zeros} entails $\mu^*=const$ almost everywhere in $\Omega$. By comparison, it also holds $\partial_t\vphi^*=0$ in $\Hu'$, for almost every $t\geq 0$. As a consequence, we infer that $$\vphi^*(t)={\vphi_\infty}$$ almost everywhere in $\Omega$, for all $t\geq 0$. Therefore, ${\vphi_\infty}$ satisfies \eqref{conv1t} for some constant $\mu_\infty\in \R$, and then ${\vphi_\infty}\in \mathcal{S}$.

    Furthermore, convergence \eqref{essenziale} implies that, up to subsequences,
    \begin{align}
    	\vphi(t+t_n)\to {\vphi_\infty} \text{ strongly in }\Hu
    	\label{ttt}
    \end{align}
    for almost any $t\in(0,\infty)$. Moreover, the convergence \eqref{AL} allows to deduce that, for any $T>0$,
    \begin{align*}
    	\sup_{t\in[0,T]}\norm{\vphi(t+t_n)-\vphi_\infty}_{\Ld}\to 0,\quad \text{ as }n\to\infty.
    \end{align*}
    Then, recalling that $f\in C([-1,1])$ and $\normmm{\vphi}\leq 1$, by Lebesgue's dominated convergence theorem we infer that, for any $T>0$,
    \begin{align}
    	\int_\Omega f(\vphi(t+t_n))\dx\to \int_\Omega f(\vphi_\infty)\dx,\quad\forall t\in[0,T],\quad \text{as }n\to \infty.\label{above}
    \end{align}
    From \eqref{ene} we infer that
    \begin{align}
    	\frac12\norm{\nabla \vphi(t_n)}_{\bLd}^2=E(\vphi(t_n))-	\int_\Omega f(\vphi(t_n))\dx\to E_\infty-\int_\Omega f(\vphi_\infty)\dx.\label{H1conv}
    \end{align}
    Additionally, for almost any $t\geq0$, owing to \eqref{ttt}, we have
    $$
    \frac12\norm{\nabla \vphi(t+t_n)}_{\bLd}^2\to \frac12\norm{\nabla \vphi_\infty}_{\bLd}^2,
    $$
    but, using \eqref{ene} and \eqref{above} once more, we also have
     $$
    \frac12\norm{\nabla \vphi(t+t_n)}_{\bLd}^2\to  E_\infty-\int_\Omega f(\vphi_\infty)\dx,\quad \forall t\geq0,
    $$
    and the uniqueness of the limit implies
    \begin{align}
    \label{constant}
    \frac12\norm{\nabla \vphi_\infty}_{\bLd}^2=E_\infty-\int_\Omega f(\vphi_\infty)\dx.
    \end{align}
    Therefore, we obtain from \eqref{H1conv} that
       \begin{align}
    	\frac12\norm{\nabla \vphi(t_n)}_{\bLd}^2\to  \frac12\norm{\nabla \vphi_\infty}_{\bLd}^2,\label{H1conv2}
    \end{align}
    which, together with the fact that $\overline{\vphi(t_n)}=\overline\vphi_0$ by mass conservation, gives by Poincaré's inequality the convergence
    \begin{align}
    \norm{ \vphi(t_n)}_{\Hu}\to  \norm{ \vphi_\infty}_{\Hu},\quad\text{ as }n\to\infty.\label{H1conv3a}
    \end{align}
    On the other hand, we know that $\vphi(t_n)\rightharpoonup\vphi_\infty$ weakly in $H^1(\Omega)$. Thus, we conclude
      \begin{align}
     \vphi(t_n)\to \vphi_\infty,\quad\text{strongly in }\Hu,\quad\text{ as }n\to\infty.\label{H1conv3}
    \end{align}
    We have thus proven that, given $\vphi_\infty\in \omega(\vphi)$ there exists a sequence $\{t_n\}_{n\in\N}$, with $t_n\to\infty$, such that \eqref{H1conv3} holds, so that the trajectories of $\vphi(\cdot)$ are precompact in $\Hu$ and $\omega(\vphi)$ can be characterized as in \eqref{omegal}.

    A straightforward consequence of the precompactness of trajectories is then property \eqref{convergence}, since, if this does not hold, by contradiction we would have that there exists $\varepsilon>0$ and a sequence $\{t_n\}_{n\in\N}$, with $t_n\to\infty$, such that
    \begin{align}
    	\inf_{\varphi_\infty\in \omega(\varphi)}\norm{\varphi(t_n)-\varphi_\infty}_{\Hu}>\varepsilon,\quad \forall n\in\N,\label{ess}
    \end{align}
    but $\{\varphi(t_n)\}_{n\in\N}$ is uniformly bounded in $H^1(\Omega)$, thus there exists a (nonrelabeled) subsequence such that $\varphi(t_n)\rightharpoonup \varphi_\infty$,  in $H^1(\Omega)$ for some $\varphi_\infty\in H^1(\Omega)$, so that $\vphi_\infty \in \omega(\vphi)$. The previous argument (see \eqref{H1conv3}) entails that there exists a subsequence such that $\vphi(t_n)\to \vphi_\infty$ strongly in $H^1(\Omega)$, which contradicts \eqref{ess}. This also immediately gives that $\omega(\vphi)$ is compact in $H^1(\Omega)$.

We now show the uniform strict separation properties of the $\omega$-limit.  It is enough to show that, given $\varphi_\infty\in \mathcal S$, there exists $\delta_{\vphi_\infty}>0$, possibly depending on $\varphi_\infty$ such that
\begin{align*}
    \norm{\vphi_\infty}_{L^\infty(\Omega)}\leq 1-\delta_{\vphi_\infty},
\end{align*}
which is trivially seen from \eqref{conv1t}, since we easily get $\norm{F'(\vphi_\infty)}_{L^\infty(\Omega)}\leq C(1+\normmm{\mu_\infty})$. Then, since $\omega(\vphi)\subset \mathcal S$, the same property holds for any element of the $\omega$-limit. To prove that the separation property is actually uniform over  $\omega(\vphi)$, we first need to show that  $\omega(\vphi)$ is bounded in $H^2(\Omega)$. By the energy inequality \eqref{energyineq} it is immediate to deduce that there exists $C>0$, only depending on the initial datum, such that
\begin{align}
\norm{\nabla \varphi_\infty}_{\bLd}\leq\liminf_{n\to\infty}\norm{\nabla \vphi(t_n)}_{\bLd} \leq \sup_{t\geq0}\norm{\nabla \vphi(t)}_{\bLd}\leq C(1+E(\vphi_0)),\quad \forall \varphi_\infty\in \omega(\vphi).\label{f}
\end{align}
Therefore, given $\varphi_\infty\in \omega(\vphi)$, by multiplying \eqref{conv1t} by $-\Delta\varphi_\infty$, integrating over $\Omega$ and after an integration by parts, we get
\begin{align*}
    \norm{\Delta \varphi_\infty}_{\Ld}^2+\int_\Omega F''(\varphi_\infty)\normmm{\nabla \varphi_\infty}^2\dx=\theta_0\int_\Omega \normmm{\nabla \varphi_\infty}^2+\int_\Omega \nabla \mu_\infty\cdot \nabla \varphi_\infty,
\end{align*}
so that, since $\mu_\infty$ is a constant and $F''\geq \theta$, we infer
$$
\norm{\Delta \varphi_\infty}_{\Ld}^2\leq \theta_0 \normm{\nabla \varphi_\infty}^2_{\bLd}\leq \theta_0C^2,\quad \forall \varphi_\infty\in \omega(\vphi),
$$
where we used \eqref{f}. Since, by mass conservation, $\overline {\vphi}_\infty=\overline{\vphi}_0$ for any $\varphi_\infty\in \omega(\vphi)$, we can deduce that $\omega(\vphi)$ is bounded in $H^2(\Omega)$. As a consequence, due to the compact embedding $H^2(\Omega)\hookrightarrow \hookrightarrow C^\alpha(\overline\Omega)$ for some $\alpha\in(0,1)$, by a simple contradiction argument we get the uniform strict separation property of $\omega(\vphi)$, i.e., \eqref{sepaglobal} (see, for instance, \cite{AW} or \cite[Proof of Lemma 3.11]{GGPS}).

The proof is concluded.
\section{Proof of Lemma \ref{twoparts}}\label{sec_prooftwoparts}
Here we assume $\Omega\subset \R^3$, the two-dimensional case being analogous. Let $\vphi$ be a global weak solution departing from the initial datum $\vphi_0\in \mathcal H_m$ (see Theorem \ref{weak}), together with its $\omega$-limit $\omega(\vphi)$ as in \eqref{omegal}. We divide the proof into two steps.

\textbf{Step 1. }
Let us first exploit the uniform strict separation of the $\omega$-limit, say \eqref{sepaglobal}. As first developed, for a completely different scope, in \cite{HKP}, we introduce the following sets, for any $\varphi_\infty\in \omega(\varphi)$,
\begin{align*}
	A_\delta(t) &:= \{x\in\Omega:\; |\varphi(x,t)|\geq 1-{\delta_1}\}, \\
	B_\delta^{\varphi_\infty}(t) &:= \{x\in\Omega:\; |\varphi(x,t)-\varphi_\infty(x)|\geq \delta_1\}.
\end{align*}
Owing to \eqref{sepaglobal}, we note that
\begin{align*}
	1- \delta_1 \leq |\varphi(x,t)| \leq |\varphi_\infty(x,t)| + |\varphi(x,t)-\varphi_\infty(x)| \leq 1-2\delta_1+ |\varphi(x,t)-\varphi_\infty(x)|
\end{align*}
for all $t\geq0$ and $x\in A_\delta(t)$. In particular, this implies
\begin{align*}
	|\varphi(x,t)-\varphi_\infty(x)|\geq \delta_1
\end{align*}
for all $t\geq 0$ and $x\in A_\delta(t)$. Thus, we get
\begin{align*}
	A_\delta(t) \subset B_\delta^{\varphi_\infty}(t) \qquad \text{ for all }t\geq 0.
\end{align*}
Using Chebyshev's inequality, we deduce
\begin{align}\label{Est:Chebyshev}
	|A_\delta(t)| \leq \int_{B_\delta^{\varphi_\infty}(t)} 1\dx \leq \int_{B_\delta^{\varphi_\infty}(t)} \frac{|\varphi(t)-\varphi_\infty|^2}{\delta_1^2}\dx \leq \frac{\norm{\varphi(t)-\varphi_\infty}^2_{\Ld}}{\delta^2_1},
\end{align}
for any $t\geq0$. Since this result holds \textit{for any }$\varphi_\infty\in \omega(\varphi_\infty)$, we can take the infimum and obtain
\begin{align}\label{Est:Chebyshev1}
	|A_\delta(t)| \leq \frac{1}{\delta^2_1}{\inf_{\varphi_\infty\in \omega(\varphi)}\norm{\varphi(t)-\varphi_\infty}^2_{\Ld}}=\frac1{\delta_1^2}\dist_{\Ld}^2(\vphi(t),\omega(\vphi))\to 0
\end{align}
as $t\to\infty$, thanks to \eqref{convergence}. Then, \eqref{conerg} holds. As a consequence, we deduce that, for any $\xi>0$, there exists $\overline T=\overline T(\xi,\delta_1)$, such that

\begin{align}\label{final}
	|A_\delta(t)| \leq \xi,\quad \forall t\geq \overline T.
\end{align}

\textbf{Step 2. De Giorgi's iterations.} We define, for $M>0$ fixed and for any $T>0$, the set of  ``good" times
$$
A_M(T):=\{t\geq T:\ \norm{\nabla \mu(t)}_{\bLd}\leq M \},
$$
which is measurable since $\nabla\mu\in L^2(0,\infty;\bLd)$.
We now perform a De Giorgi's iteration scheme to prove the validity of the strict separation property on $A_M(T)$, for $T>0$ sufficiently large. The proof takes inspiration from \cite{GalP}. Let us then fix $0<\delta\leq \delta_1$, for $\delta_1>0$ given in \eqref{sepaglobal}. We define, as usual in this argument, the sequence
\begin{align}
	k_n=1-\delta-\frac{\delta}{2^n}, \quad \forall n\geq 0,
	\label{kn}
\end{align}
where
\begin{align}
	1-2\delta< k_n<k_{n+1}<1-\delta,\qquad \forall n\geq 1,\qquad k_n\to 1-\delta\qquad \text{as }n\to \infty.
	\label{kn1}
\end{align}
We then set
\begin{align}
	\varphi_n(x,t):=(\varphi-k_n)^+.
	\label{phik0}
\end{align}
In conclusion, we define
$$
y_n(t)=\int_{A_n(t)}1dx,\qquad \forall n\geq0,\quad \forall t\geq0,
$$
where
$$
A_n(t):=\{x\in \Omega: \vphi(x,t)\geq k_n\},\quad\forall t\geq0.
$$
Now we fix $\xi$ (and thus $\overline T(\xi,\delta_1)>0$) such that (see \eqref{final})
$$
\normmm{A_\delta(t)}\leq \frac{\normmm{\Omega}}{8},\quad \quad \forall t\geq \overline T.
$$
Then, recalling that $0\leq \varphi_n\leq 2\delta$ (see \cite{P}), with $\delta\in(0,\delta_1)$, we find
$$
\overline\varphi_n(t)=\frac{\int_{A_n(t)} \vphi_n(t)\dx}{\normmm{\Omega}}\leq 2\delta\frac{\normmm{A_n(t)}}{\normmm{\Omega}}\leq 2\delta\frac{\normmm {A_\delta(t)}}{\normmm{\Omega}}\leq  \frac{\delta}4,\quad \forall t\geq \overline T.
$$
Therefore, we have
\begin{align}
	0<k_n+\overline\varphi_n(t) \leq 1-\frac{3\delta}4-\frac{\delta}{2^n},\quad \forall t\geq \overline T,
	\label{ext}
\end{align}
so that the quantity $F'(k_n+\overline\varphi_n(t))\leq F'(1-\frac{3\delta}4)<+\infty$ used below is well defined and finite for any $t\geq \overline T$, since $\delta<\delta_1$ is fixed.

Now, let us fix $t\geq \overline T$. For any $n\geq 0$, we consider the test function $v =
\vphi_{n}-\overline\vphi_n$, multiply equation \eqref{A1} by $v$, and integrate over $\Omega$.
After an integration by parts, taking into account the boundary conditions,
we obtain:
\begin{equation}
	\Vert \nabla \vphi_{n}\Vert^{2}_{\Ld} + \int_{\Omega} F^{\prime}(\vphi) (\vphi_{n}-\overline\vphi_n)dx =
	\theta_{0} \int_{\Omega} \vphi (\vphi_{n}-\overline\vphi_n)dx + \int_{\Omega} \mu (\vphi_{n}-\overline\vphi_n)dx,\label{base}
\end{equation}
for any $t\in [\overline T,\infty)$. Here, we used the identity:
\begin{equation}
	\int_{A_{n}} \nabla \vphi \cdot \nabla \vphi_{n}dx = \Vert \nabla
	\vphi_{n}\Vert^{2}_{\bLd}.  \label{Fss}
\end{equation}
Observe now that
\begin{align*}
	&\int_{\Omega}F^\prime(\varphi(t))(\varphi_n(t)-\overline{\varphi}_n(t))dx\\&
	=\int_{\{x:\ \varphi_n(x,t)= 0\}}F^\prime(\varphi(t))(\varphi_n(t)-\overline{\varphi}_n(t))dx+\int_{\{x:\ \varphi_n(x,t)>0\}}F^\prime(\varphi(t))(\varphi_n(t)-\overline{\varphi}_n(t))dx\\&:=
	J_1+J_2.
\end{align*}
Concerning $J_1$, for any $t\geq 0$, we have
$$
Z_n(t)=\{x\in\Omega:\ \varphi_n(x,t)= 0\}=\{x\in\Omega:\ \varphi(x,t)\leq k_n\},
$$
and for any $x\in Z_n(t)$, since $F'$ is monotone increasing, with $F'(s)\leq 0$ for $s\leq 0$, there holds
\begin{align*}
	&{ F^\prime(\varphi)(\varphi_n-\overline{\varphi}_n)}=-{F^\prime(\varphi)\overline{\varphi}_n}\\&
	= \underbrace{-{F^\prime(\varphi)\overline{\varphi}_n}\chi_{\{\varphi\leq 0\}}}_{\geq0}-{F^\prime(\varphi)\overline{\varphi}_n}\chi_{\{0<\varphi\leq k_n\}}\\&\geq -{F^\prime(\varphi)\overline{\varphi}_n}\chi_{\{0<\varphi\leq k_n\}}\geq -{F^\prime(k_n)\overline{\varphi}_n}\chi_{\{0<\varphi\leq k_n\}}.
\end{align*}
As a consequence, we can write
\begin{align}
	\nonumber J_1(t)&=\int_{\{x:\ \varphi_n(x,t)= 0\}}F^\prime(\varphi)(\varphi_n-\overline{\varphi}_n)dx\\&\nonumber\geq -\int_{\{x:\ \varphi_n(x,t)= 0\}}{F^\prime(k_n)\overline{\varphi}_n}\chi_{\{0<\varphi\leq k_n\}}dx\\&
	=-\int_{\{x:\ 0<\varphi(x,t)\leq k_n\}}{F^\prime(k_n)\overline{\varphi}_n}dx.\label{phin}
\end{align}
Let us consider $J_2$. Recalling the definition of $\varphi_n$, we have that, for any $t\geq0$,
\begin{align*}
	\{x\in\Omega:\ \varphi_n(x,t)>0\}=\{x\in\Omega:\ \varphi(x,t)>k_n\},
\end{align*}
and thus
\begin{align*}
	J_2(t)&=\int_{\{x:\ k_n<\varphi(x,t)\leq \overline \varphi_n(t)+k_n\}}F^\prime(\varphi)(\varphi_n-\overline{\varphi}_n)dx\\&\quad +\int_{\{x:\ \varphi(x,t)>\overline \varphi_n(t)+k_n\}}F^\prime(\varphi)(\varphi_n-\overline{\varphi}_n)dx=:J_3(t)+J_4(t).
\end{align*}
Note that $J_3$ can be treated as follows. Recalling the monotonicity of $F'$, that $t\geq \overline T$, and \eqref{ext}, we deduce
\begin{align*}
	J_3(t)&=\int_{\{x:\ k_n<\varphi(x,t)\leq \overline \varphi_n(t)+k_n\}}F^\prime(\varphi)(\varphi-k_n-\overline{\varphi}_n)dx\\&=-\int_{\{x:\ k_n<\varphi(x,t)\leq \overline \varphi_n(t)+k_n\}}F^\prime(\varphi)\underbrace{(-\varphi+k_n+\overline{\varphi}_n)}_{\geq 0}dx\\&\geq -\int_{\{x:\ k_n<\varphi(x,t)\leq \overline \varphi_n(t)+k_n\}}F'(\overline\varphi_n+k_n)(-\varphi+k_n+\overline{\varphi}_n)dx.
\end{align*}
Analogously, we have
\begin{align*}
	J_4(t)&=\int_{\{x:\ \varphi(x,t)>\overline \varphi_n(t)+k_n\}}F^\prime(\varphi)\underbrace{(\varphi-k_n-\overline{\varphi}_n)}_{\geq0 }dx\\&\geq \int_{\{x:\ \varphi(x,t)>\overline \varphi_n(t)+k_n\}}F^\prime(\overline \varphi_n+k_n){(\varphi-k_n-\overline{\varphi}_n)}dx.
\end{align*}
Therefore, summing up the two estimates, we get
\begin{align*}
	J_3(t)+J_4(t)&\geq \int_{\{x:\ k_n<\varphi(x,t)\leq \overline \varphi_n(t)+k_n\}}F'(\overline\varphi_n+k_n)(\varphi-k_n-\overline{\varphi}_n)dx\\&\quad +\int_{\{x:\ \varphi(x,t)>\overline \varphi_n(t)+k_n\}}F^\prime(\overline \varphi_n+k_n){(\varphi-k_n-\overline{\varphi}_n)}dx\\&
	=\int_{\{x:\ k_n<\varphi(x,t)\}}F'(\overline\varphi_n+k_n)(\varphi-k_n-\overline{\varphi}_n)dx.
\end{align*}
If we set $\{x\in\Omega:\ k_n<\varphi(x,t)\}=W_n(t)$, then we infer
\begin{align*}
	&\int_{\{x:\ k_n<\varphi(x,t)\}}F'(\overline\varphi_n+k_n)(\varphi-k_n-\overline{\varphi}_n)dx\\&=F'(\overline\varphi_n+k_n)\int_{W_n(t)}(\varphi_n-\overline{\varphi}_n)dx\\&\quad \pm
	F'(\overline\varphi_n+k_n)\int_{\Omega\setminus W_n(t)}(\varphi_n-\overline{\varphi}_n)dx\\&=\underbrace{F'(\overline\varphi_n+k_n)\int_{\Omega}(\varphi_n-\overline{\varphi}_n)dx}_{=0}-F'(\overline\varphi_n+k_n)\int_{\{x:\ \varphi(x,t)\leq k_n\}}\overline{\varphi}_ndx.
\end{align*}
To sum up, collecting the estimates of $J_1$-$J_4$, we have obtained that
\begin{align}
	&\nonumber\int_{\Omega}F^\prime(\varphi)(\varphi_n-\overline{\varphi}_n)dx\\&\geq -\int_{\{x:\ 0<\varphi(x,t)\leq k_n\}}{F^\prime(k_n)\overline{\varphi}_n}dx-F'(\overline\varphi_n+k_n)\int_{\{x:\ \varphi(x,t)\leq k_n\}}\overline{\varphi}_ndx,\quad \forall t\geq \overline T.
	\label{controls}
\end{align}
Using now Poincaré's (with constant $C_P>0$) and H\"older's inequalities, and recalling that $0\leq
\vphi_n\leq2\delta$, we get
\begin{align}
	& \normmm{\int_{\Omega} \mu (\vphi_{n}-\overline\vphi_n)dx}=\normmm{\int_{\Omega} (\mu-\overline\mu) \vphi_{n}dx}  \leq \Vert \vphi_{n}\Vert_{L^{\infty}(\Omega)} \Vert \mu-\overline\mu
	\Vert_{L^{6}(\Omega)} \left(\int_{A_{n}} 1dx\right)^{\frac56}  \notag
	\\
	& \leq 2\delta C_P\Vert \nabla \mu \Vert y_{n}^{\frac 56}.  \label{c2}
\end{align}
Here we need our definition of ``good" times. Indeed, for any $t\in A_M(T)$, for $T>\overline T$, it holds
\begin{align}
	& \normmm{\int_{\Omega} \mu (\vphi_{n}-\overline\vphi_n)dx}     \notag
	 \leq 2\delta C_P\Vert \nabla \mu \Vert y_{n}^{\frac 56}\leq 2C_PM\delta y_{n}^{\frac 56},\quad \forall t\in A_M(T).  \label{c3}
\end{align}
Using the above results above, from \eqref{base} we infer
\begin{align}
\nonumber	\norm{\nabla \vphi_n}^2_{\bLd}\nonumber&\leq 2C_PM\delta y_{n}^{\frac 56}+	\theta_{0} \int_{\Omega} \vphi (\vphi_{n}-\overline\vphi_n)dx+\left({F^\prime(k_n)}+F'(\overline\varphi_n+k_n)\right)\int_{\Omega}\overline{\varphi}_ndx\\&
\leq C_A\delta\left(1+F'(1-\delta)+F'(1-\frac34\delta)\right)y_n^{\frac56},
\end{align}
where we recall that
\begin{align*}
	\int_\Omega\overline\vphi_n(t)\dx=\int_\Omega\frac{\int_{A_n(t)}\vphi_n(y,t)\d y}{\normmm{\Omega}}\dx\leq {2\delta}y_n(t)\leq 2\delta\normmm{\Omega}^{\frac16}y_n(t)^\frac56,
\end{align*}
and ($\normmm{\vphi}< 1$ almost everywhere in $\Omega\times(0,\infty)$),
\begin{align*}
    \normmm{\theta_{0} \int_{\Omega} \vphi (\vphi_{n}-\overline\vphi_n)dx}\leq 2\theta_0\int_\Omega \vphi_n\dx\leq 4\delta\theta_0 y_n(t)\leq  4\delta\theta_0\normmm{\Omega}^{\frac16}y_n(t)^\frac56.
\end{align*}
Then, we can write
\begin{align}
	\norm{\nabla \vphi_n(t)}^2_{\bLd}\leq C_B(\delta,M)y_n^{\frac56}(t),\quad \forall t\in A_M(T),
\end{align}
for some $C_B>0$, depending only on $\delta$ and $M$.
On the other hand, for any $t\geq 0 $ and for almost any $x\in A_{n+1}(t)$, we have
\begin{align}
	&\nonumber\varphi_n(x,t)=\varphi(x,t)-\left[1-\delta-\frac{\delta}{2^n}\right]\\&=
	\underbrace{\varphi(x,t)-\left[1-\delta-\frac{\delta}{2^{n+1}}\right]}_{\varphi_{n+1}(x,t)\geq 0}+\delta\left[\frac{1}{2^{n}}-\frac{1}{2^{n+1}}\right]\geq \frac{\delta}{2^{n+1}},\label{basic}
\end{align}
which implies
\begin{align*}
	\int_{\Omega}\vert\varphi_n\vert^3dx\geq \int_{A_{n+1}(t)}\vert\varphi_n\vert^3dx\geq \left(\frac{\delta}{2^{n+1}}\right)^3\int_{A_{n+1}(t)}dx=\left(\frac{\delta}{2^{n+1}}\right)^3y_{n+1}.
\end{align*}
Then we have
\begin{align}
	&\nonumber\left(\frac{\delta}{2^{n+1}}\right)^3y_{n+1}\leq \int_{\Omega}\vert\varphi_n\vert^3dx\\&=\int_{A_n(t)}\vert\varphi_n\vert^3dx \leq \left(\int_{\Omega}\vert\varphi_n\vert^{\frac{10}{3}}dx\right)^{\frac{9}{10}}\left(\int_{A_n(t)}1\ dx\right)^{\frac{1}{10}}.
	\label{est2}
\end{align}
Notice that, by Gagliardo-Nirenberg's inequalities, we get
\begin{align*}
	&\int_{\Omega}\vert\varphi_n\vert^{\frac{10}{3}}dx\leq \hat{C} \Vert\varphi_n\Vert_{H^1(\Omega)}^{2}\Vert\varphi_n\Vert^{\frac{4}{3}}_{\Ld}\leq \hat{C} \left(\Vert\varphi_n\Vert^2_{\Ld}+\Vert\nabla\varphi_n\Vert^2_{\Ld}\right)\Vert\varphi_n\Vert^{\frac{4}{3}}_{\Ld}.
\end{align*}
On account of \eqref{est} and recalling that $\normmm{\vphi_n}\leq 2\delta$ and $y_n\leq \normmm{\Omega}$, for any $t\in A_M(T)$, $T>\overline T$,
we find
\begin{align*}
	&\int_{\Omega}\vert\varphi_n\vert^{\frac{10}{3}}dx\leq \hat{C}\Vert\varphi_n(t)\Vert_{\Ld}^{\frac{4}{3}}\left(\Vert\varphi_n(t)\Vert^2_{\Ld}+\Vert\nabla\varphi_n(t)\Vert^2_{\bLd}\right)\\&
\leq   \hat{C}(2\delta)^\frac43y_n^\frac23((2\delta)^2y_n+C_B(\delta,M)y_n^{\frac56})
\\&
	\leq C_C(\delta,M)y_n^{\frac{3}{2}},
\end{align*}
for some $ C_C(\delta,M)>0$, where we used
\begin{align*}
    \Vert\varphi_n(t)\Vert_{\Ld}=\left(\int_{A_n(t)}\vphi_n^2(x,t)\dx\right)^\frac12\leq 2\delta y_n(t)^\frac12.
\end{align*}
Therefore, we infer from \eqref{est2} that
\begin{align}
	&\nonumber\left(\frac{\delta}{2^{n+1}}\right)^3y_{n+1}\leq \left(\int_{\Omega}\vert\varphi_n\vert^{\frac{10}{3}}dx\right)^{\frac{9}{10}}\left(\int_{A_n(t)}1\ dx\right)^{\frac{1}{10}}\\&\leq C_C(\delta,M)^\frac9{10}y_n^{\frac{29}{20}},
	\label{est3}
\end{align}
for any $t\in A_M(T)$. In conclusion, we end up with
\begin{align}
	y_{n+1}(t)\leq \frac1{\delta^3}2^{{3n+3}}C_C(\delta,M)^\frac{9}{10}y_n^{\frac{29}{20}}(t),\qquad \forall n\geq 0,
	\label{last0}
\end{align}
for any $t\in A_M(T)$.
Thus we can apply the well known geometric Lemma \ref{conv}. In particular, using the notation of the lemma, we have $b=2^3>1$, $C=\frac1{\delta^3}2^{{3}}C_C(\delta,M)^\frac{9}{10}>0$, $\varepsilon=\frac{9}{20}$, to get that ${y}_n\to 0$, as long as
$$
{y}_0(t)\leq C^{-\frac{20}{9}}b^{-\frac{400}{81}},
$$
i.e.,
\begin{align}
	y_0(t)\leq C_D(\delta,M),
	\label{last}
\end{align}
for some $C_D(\delta,M)>0$ which is \textit{independent of} $t$.  Then, from \eqref{final}, for any $\xi>0$ there exists $\overline T(\xi)$, possibly larger, so that, since $\delta\in(0,\delta_1)$, we have
\begin{align*}
	&y_0(t)=\int_{A_0(t)}1dx\leq\int_{\{x\in\Omega:\ \varphi(x,t) \geq 1-2\delta\}}1dx\leq \normmm{A_\delta(t)}\leq \xi,\quad \forall t\geq \overline T.
\end{align*}
Therefore, if we choose $\xi$ sufficiently small (and thus we fix $\overline T(\xi)$) such that
$$
\xi\leq C_D(\delta,M),
$$
then \eqref{last} holds for any $t\in A_M(T)$, with $T>\overline T$.
In the end, passing to the limit in $y_n(t)$ as $n\to\infty$, we have obtained that
$$
\Vert(\varphi(t)-(1-\delta))^+\Vert_{L^\infty(\Omega)}=0,
$$
for any $t\in A_M(T)$, with $T>\overline T$. Since $\delta\in(0,\delta_1)$ does not depend on $t$, we then have
$$
\sup_{t\in A_M(T)}\Vert(\varphi(t)-(1-\delta))^+\Vert_{L^\infty(\Omega)}=0.
$$
We now repeat the very same argument for the case $(\varphi-(-1+\delta))^-$ (using $\varphi_n(t)=(\varphi(t)+k_n)^-$) to get
$$
\sup_{t\in A_M(T)}\Vert(\varphi(t)-(-1+\delta))^-\Vert_{L^\infty(\Omega)}=0.
$$
Therefore, we have shown that there exist $\delta\in(0,\delta_1)$ and $T_S>\overline T>0$ such that
\begin{align}
\sup_{t\in A_M(T_S)}\norm{\varphi(t)}_{L^\infty(\Omega)}\leq 1-\delta.
	\label{end}
\end{align}
The proof of \eqref{asympt} is thus concluded.

\section{Proof of Theorem \ref{uniqueeq}}\label{sec:proofLoja}
We first recall the \L ojasiewicz-Simon inequality we need (see \cite[Proposition 6.1]{AW}), which is valid thanks to the regularity $C^2$ of the boundary $\partial\Omega$:
	\begin{proposition}
				 Assume that $F$ is additionally real analytic in $(-1,1)$ {and let} $\vphi\in  \mathcal H_k$ {be} such that $-1+\gamma\leq \vphi(x)\leq 1-\gamma$, for {almost }any $x\in \Omega$ {and }for some $\gamma\in(0,1)$. {Furthermore, let } $\vphi_\infty\in \mathcal{S}$ {be fixed} such that $-1+\gamma\leq \vphi_\infty(x)\leq 1-\gamma$ for any $x\in\Omega$. Then there exist $\vartheta\in \left(0,\frac{1}{2}\right)$, $\eta>0$ and a positive constant $C$ such that
				\begin{align}
					\vert{{E}}(\vphi)-{E}(\vphi_\infty)\vert^{1-\vartheta}\leq C\| \delta E(\varphi)\|_{\Hu'},
					\label{ener}
				\end{align}
				{provided that } $\|\vphi-\vphi_\infty\|_{\Hu}\leq \eta$,
				\label{Lojaw}
where $\delta E:H^1_{(k)}(\Omega)\to H^1_{(0)}(\Omega)'$ is the Frechét derivative of $E:H^1_{(k)}(\Omega)\to \R$.
\end{proposition}
Thanks to \eqref{energyineq} and \eqref{ene}, we have that ${{E}}(\vphi(t))\geq {E}_\infty$ and that ${E}(\vphi(t))\to {E}_\infty$, as $t\to\infty$. Let us now fix $M>0$. Then we choose $\gamma$ equal to the value of $\delta$ given in Lemma \ref{twoparts} (see inequality \eqref{asympt}), so that, as $\gamma\leq\delta_1$, it holds, by \eqref{sepaglobal}, $-1+\gamma\leq \vphi_\infty\leq 1-\gamma$ in $\Omega$, for any $\vphi_\infty\in\omega(\vphi)$. Moreover, for any $\vphi_{\infty,m}\in \omega(\vphi)$ we can find $\theta_m\in \left(0,\frac{1}{2}\right)$ and $\eta_m>0$, given {by }Proposition \ref{Lojaw}, for which \eqref{Lojaw} is valid with a constant $C_m$. From Lemma \ref{convaaa} we observe that $\omega(\vphi)$ is compact in $H^1(\Omega)$ and bounded in $H^2(\Omega)$. {We can thus find a finite family of $H^1(\Omega)$-open balls, say $\{B_{\eta_m}\}_{m=1}^{M_1}$, centered at points $\{\vphi_{\infty,m}\}_{m=1}^{M_1}\subset \omega(\vphi)$ and with radii $\eta_m$ (depending on the center $\vphi_{m,\infty}\in\omega(\vphi)$), such that
				\begin{equation*}
					\bigcup_{ \varphi_\infty \in \omega(\varphi)} \{\varphi_\infty\} \subset U:= \bigcup_{m=1}^{M_{1}}B_{\eta_m}.
				\end{equation*}}
	Recalling \eqref{constant}, which is valid for any $\vphi_\infty\in \omega(\vphi)$, we infer that the energy functional $E(\cdot)$ is constant over $\omega(\vphi)$. Additionally, since {the centers  $\{\vphi_m\}_{m = 1}^{M_1}$ are in finite number, }we can infer that \eqref{ener} holds \textit{uniformly}, with suitable constants, for any $\vphi\in U$ such that $\|\vphi\|_{L^\infty(\Omega)}\leq 1-\gamma$, and we can substitute $E(\vphi_\infty)$ with $E_\infty$. 
    
    By \eqref{convergence}, we deduce that there exists ${t_*}>0$ such that $\vphi(t)\in U$ for any $t\geq t_*$. Additionally, observe that, thanks to \eqref{asympt}, there exists $T_S>0$ such that the uniform strict separation property holds on the set of ``good" times $A_M(T_S)=\{t\geq T_S:\ \norm{\nabla\mu(t)}_{\bLd}\leq M\}$, namely
                \begin{align*}
                    \sup_{t\in A_M(T_S)}\norm{\vphi(t)}_{L^\infty(\Omega)}\leq 1-\delta.
                \end{align*}
                Therefore, thanks to the choice of $\gamma=\delta$, we get from \eqref{asympt} and \eqref{ener} that
					\begin{align}
						\label{pp}\left({{E}}_{CH}(\vphi(t))-{{E}}_\infty\right)^{1-\vartheta}\leq C\norm{\mu(t)-\overline\mu(t)}_{H^1(\Omega)'} \leq C\| \nabla\mu(t)\|_{\bLd},\quad \forall t\in A_M(T_S)\cap [t_*,\infty).
					\end{align}
The core of the novel argument follows. We have, by the energy inequality \eqref{energyineq} and recalling $m(\cdot)\geq m_*>0$ since the mobility is non-degenerate,
					\begin{align*}
						&\nonumber m_*\int_s^t\norm{\nabla \mu(\tau)}_{\bLd}^2\d\tau \leq \int_s^t \int_\Omega m(\vphi(x,\tau))\normmm{\nabla\mu(x,\tau)}^2\dx\dtau
						\leq E(\vphi(s))-E(\vphi(t)),
					\end{align*}
					for any $t>0$ and almost any $s\in[0,t]$, $s=0$ included.
					This gives
					\begin{align}
						& \left(\int_s^t \norm{\nabla\mu(\tau)}_{\Ld}^2\dtau\right)^{2(1-\vartheta)}
						\leq \frac1{m_*^{2(1-\vartheta)}}(E(\vphi(s))-E(\vphi(t)))^{2(1-\vartheta)}.\label{cv}
					\end{align}
					We now let $t\to\infty$ in \eqref{cv}, and obtain, recalling that $E(\vphi(t))\to E_\infty$ as $t\to \infty$, for almost any $s\in(t_*,\infty)$,
					\begin{align}
						&\nonumber \left(\int_s^\infty \norm{\nabla\mu(\tau)}_{\Ld}^2\dtau\right)^{2(1-\vartheta)}\\&
						\leq \frac1{m_*^{2(1-\vartheta)}}(E(\vphi(s))-E_\infty)^{2(1-\vartheta)}(\chi_{A_M(T_S)}(s)+\chi_{(t_*,\infty)\setminus A_M(T_S)}(s)).\label{cv22}
					\end{align}
			Observe that, for almost any $s\in (t_*,\infty)\setminus A_M(T_S)$ (i.e., the ``bad" times), it holds $\norm{\nabla \mu(s)}_{\bLd}\geq M$. Thus, recalling that $E(\vphi(t))\leq E(\vphi_0)$ for any $t\geq0$ and $t\mapsto E(\vphi(t))$ is monotone decreasing, we infer
					$$
				(E(\vphi(s))-E_\infty)^{2(1-\vartheta)}\chi_{(t_*,\infty)\setminus A_M(T_S)}(s)\leq (2E(\vphi_0))^{2(1-\vartheta)}\frac{\norm{\nabla\mu(s)}_{\bLd}^2}{M^2}\chi_{(t_*,\infty)\setminus A_M(T_S)}(s),
				$$
				for almost any $s\in (t_*,\infty)\setminus A_M(T_S)$.
				On the other hand, on the ``good" times, thanks to \eqref{pp}, we have
					$$
				(E(\vphi(s))-E_\infty)^{2(1-\vartheta)}\chi_{ A_M(T_S)}(s)\leq C^2{\norm{\nabla\mu(s)}_{\bLd}^2}\chi_{A_M(T_S)}(s),
				$$
				for almost any $s\in  A_M(T_S)$.
				As a consequence, we deduce from \eqref{cv22} that
					\begin{align}
					& \left(\int_s^\infty \norm{\nabla\mu(\tau)}_{\bLd}^2\dtau\right)^{2(1-\vartheta)}
					\leq \frac1{m_*^{2(1-\vartheta)}}\left(C^2+\frac{(2 E(\vphi_0))^{2(1-\vartheta)}}{M^2}\right)\norm{\nabla \mu (s)}^2_{\bLd},\label{cv2}
				\end{align}
				for almost any  $s\in (t_*,\infty)$.
			    Using now Lemma \ref{Feireisl} with its notation, namely,
$$
Z(\cdot)=\norm{\nabla\mu(\cdot)}_{\bLd}, \quad \alpha=2(1-\vartheta)\in(1,2), \quad \zeta=\tfrac1{m_*^{2(1-\vartheta)}}(C^2+\tfrac{(2 E(\vphi_0))^{2(1-\vartheta)}}{M^2})>0,
$$ and $\mathcal M=(t_*,\infty)$, \eqref{fin} yields
					\begin{align}
						\nabla \mu\in L^1(t^*,\infty;\mathbf L^2(\Omega)) . \label{bA}
					\end{align}
					Thus, by comparison, we deduce that $\partial_t\vphi\in L^1(t^*,\infty;\Hu')$. Hence, we have
					$$
					\vphi(t)=\vphi(t^*)+ \int_{t^*}^t\partial_t\vphi(\tau) \: \d\tau \to {\vphi_\infty}\quad \text{ in }\Hu', \text{ as } t \to \infty,$$for some ${\vphi_\infty}\in \Hu'$. This means that $\vphi(t)$ converges in $\Hu'$ as $t\to\infty$ and, by uniqueness, we conclude that $\omega(\vphi)$ is a singleton. Using the interpolation inequality $\norm{\cdot}_{\Ld}\leq C\norm{\cdot}_{\Hu'}^\frac12\norm{\cdot}_{\Hu}^\frac12$ together with $\vphi\in L^\infty(0,\infty;\Hu)$, we find
                    $$
                    \norm{\vphi(t)-\vphi_\infty}_{\Ld}\to0,
                    $$
                    as $t\to\infty$. In order to show \eqref{equil}, it is enough to use \eqref{convergence}, as now $\omega(\vphi)$ is a singleton. The proof is complete.
	\section{The Abels-Garcke-Gr\"un system with non-degenerate mobility}
\label{proofAGG}
\subsection{Proof of Lemma \ref{lemmauconv}}\label{secproofconvu}
The proof can be obtained by following the same arguments as in \cite[Lemma 3.2]{AGGG}, which are only based on the energy inequalities
 \eqref{energyineq3} and \eqref{energyineq2}. Indeed, first observe that \eqref{energyineq2} entails $\nabla\mu\in L^2(0,\infty;\bLd)$.
  Thus, for any $\varepsilon>0$, there exists $T_1(\varepsilon)>0$ such that $\norm{\nabla\mu}_{L^2(T_1,\infty;\bLd)}\leq \varepsilon$. Let us then fix such $\varepsilon>0$. By Korn's inequality and using \eqref{energyineq2} once more, there exists $T_2(\varepsilon)$ such that $\norm{\bu(T_2)}_{\bLds}\leq \varepsilon$. As a consequence, by choosing $T>\max\{T_1,T_2\}$, recalling that $\normmm{\vphi}< 1$ almost everywhere in $\Omega\times(0,\infty)$, and using standard inequalities, we infer from \eqref{energyineq3} the following inequality
\begin{align*}
&\nonumber\frac12\int_\Omega \rho(\vphi(t))\normmm{\bu(t)}^2\dx+\int_{T}^t\int_\Omega\nu(\vphi(\tau)\normmm{D\bu(\tau)}^2\dx\d\tau\\&\leq \frac12\int_\Omega \rho(\vphi(T))\normmm{\bu(T_2)}^2\dx-\int_T^t\int_\Omega \bu(\tau)\cdot\nabla\mu(\tau)\vphi(\tau)\dx\d\tau\\&
\leq C_1\varepsilon+\norm{\bu}_{L^2(T,\infty;\bLds)}\norm{\nabla\mu}_{L^2(T,\infty;\bLd)}\leq (C_1+C_2)\varepsilon,\quad \forall t\geq T(\varepsilon),
    	\end{align*}
        with $C_1,C_2>0$ independent of $T$ and $\varepsilon$.
This gives that $\norm{\bu(t)}_{\bLds}\to 0$ as $t\to\infty$.
\subsection{Proof of Theorem \ref{uniqueeq1}}
Let us state and prove the following preliminary lemma.
\label{secprooffinal}
    \begin{lemma}
				\label{convaaa1}
				Let the assumptions of Theorem \ref{weakAGG} hold. We have
				$$
				    \omega(\bu,\vphi)\subset \mathcal{S}_1.
				$$
				Moreover, $\omega(\bu,\vphi)$ is bounded in $\{ \mathbf 0\} \times H^2(\Omega)$, and there exists $\delta_1>0$ such that
				\begin{align}
				\| \vphi_\infty\|_{L^\infty(\Omega)}\leq 1-2\delta_1,\quad \forall \: \vphi_\infty\in \omega(\vphi).
				\label{sepaglobal12
                }
                \end{align}
The trajectories of $\vphi(\cdot)$ are precompact in $H^1(\Omega)$. Moreover, we have that $\omega(\vphi)$ is compact in $H^1(\Omega)$, it holds  
                \begin{align}
\nonumber&\omega(\bu,\varphi)\\&=\{(\mathbf 0,\widetilde{\varphi})\in  \{\mathbf 0\}\times \mathcal{H}_k:\exists t_n\to \infty \text{ s.t. }\varphi(t_n)\to \widetilde{\varphi}\text{ in }H^1(\Omega)\text{ and }\bu(t_n)\to\mathbf 0\text{ in }\bLds\},\label{omegal12}
\end{align}
and 
\begin{align}
	\lim_{t\to\infty}\dist_{ \bLds\times H^1(\Omega)}((\bu(t),\varphi(t)),\omega(\bu,\vphi))=0.\label{convergence2}
\end{align}
			\end{lemma}
\begin{proof}
    The proof of the lemma follows closely the one of the analogous Lemma \ref{convaaa}, up to the additional presence of the advective velocity $\bu$ in the Cahn--Hilliard equation. We thus only highlight the adaptations.
   Let us consider a sequence $t_n\to \infty$ such that $\varphi(t_n)\rightharpoonup  \widetilde{\varphi}$ weakly in $H^1(\Omega)$, with $\widetilde\vphi\in \omega(\varphi)$ (and thus, up to subsequences, $\vphi(t_n)\to\widetilde\vphi$ strongly in $\Ld$). We then consider the sequence of trajectories $\bu_n(t)=\bu(t+t_n)$, $\varphi_n(t):=\varphi(t+t_n)$, $\mu_n(t):=\mu(t+t_n)$ and we observe that
   \begin{align}
      \label{eqs1b}  &\langle \partial_t\vphi_n,v\rangle-(\bu_n\vphi_n,\nabla v)+(m(\vphi_n)\nabla \mu_n,\nabla v)=0,\quad \forall v\in H^1(\Omega), \text{ for a.a. }t\geq 0,\\&
        \mu_n=-\Delta\vphi_n+f'(\vphi_n)\quad\text{ a.e. in }\Omega\times(0,\infty).\label{eqs2b}
    \end{align}

   By the energy inequality \eqref{energyineq2} we get that ${{E}}_{tot}(\bu(t_n),\vphi(t_n))\leq {{E}}_{tot}(\bu_0,\vphi_0)$ for any $n$. As a consequence, for any $T>0$, arguing as for \eqref{energy}, we can find a constant $C(T)>0$ independent of $n$ such that
				\begin{align}
				\label{energia}\norm{\bu_n}_{L^\infty(0,T;\bLds)\cap L^2(0,T;\mathbf H^1(\Omega))}+\| \vphi_n\|_{L^\infty(0,T;H^1(\Omega))}+\| \mu_n\|_{L^2(0,T;H^1(\Omega))}\leq C(T).
				\end{align}
Then, following line by line the proof of \eqref{elliptic}, we get
\begin{align}
\norm{F'(\vphi_n)}_{L^2(0,T;L^2(\Omega))}+\norm{\vphi_n}_{L^2(0,T;H^2(\Omega))}\leq C(T).
    \label{elliptic2}
\end{align}
Recalling \eqref{energia} and $\normmm{\vphi_n}\leq 1$ almost everywhere in $\Omega\times(0,\infty)$, we see that
$$
\normmm{\int_\Omega \vphi_n\bu_n\cdot   \nabla v}\leq \norm{\vphi_n}_{L^\infty(\Omega)}\norm{\bu_n}_{\bLds}\normm{v}_{H^1(\Omega)}\leq \norm{\bu_n}_{\bLds}\normm{v}_{H^1(\Omega)},\quad \forall v\in H^1(\Omega).
$$
Thus, by comparison in \eqref{eqs1b}, using also the bounds \eqref{energia}, we infer
\begin{align}
\norm{\partial_t\vphi_n}_{L^2(0,T;H^1(\Omega)')}\leq C(T).
    \label{timeder}
\end{align}
As a consequence of \eqref{uconv} and \eqref{energia}-\eqref{timeder}, we deduce that all the convergences \eqref{est}-\eqref{AL} also hold also in the present case. Therefore, letting $n\to\infty$, we  infer
that {the limit pair } $(\vphi^*,\mu^*)$ satisfies, for any $T>0$,
				\begin{align}
					\label{c1}
&\langle \partial_t\vphi^*,v\rangle+(m(\vphi^*)\nabla \mu^*, \nabla v)=0,\quad \forall \: v\in \Hu,\quad \text{a.e. in }(0,T),\\
					&\mu^*=-\Delta\vphi^*+f'(\vphi^*),\quad\text{a.e. in }\Omega\times(0,T),\\
&\partial_{\bn}\varphi_\infty=0, \quad\text{a.e. on } \partial\Omega\times(0,T),
				\end{align}
				with initial datum $\vphi^*(0)={\vphi_\infty}$. Then, since $E_{tot}(\bu(\cdot),\vphi(\cdot))$ is nonincreasing in time, there exists $E_\infty$ so that $E_{tot}(\bu(t),\vphi(t))\to E_\infty$ as $t\to\infty$. By the same argument used to get \eqref{zeros}, we deduce that $\mu_\infty=const$, which by comparison in \eqref{c1} gives $\vphi^*$ constant in time, coinciding with $\vphi_\infty$. The remaining part of the proof then goes as in the proof of Lemma \ref{convaaa}, by replacing \eqref{H1conv} with
                \begin{align*}
    	&\frac12\norm{\nabla \vphi(t_n)}_{\bLd}^2\\&=E_{tot}(\bu(t_n),\vphi(t_n))-\frac12\int_\Omega\rho(\vphi(t_n))\normmm{\bu(t_n)}^2\dx-	\int_\Omega f(\vphi(t_n))\dx\to E_\infty-\int_\Omega f(\vphi_\infty)\dx,
    \end{align*}
        as $n\to\infty,$ since (see \eqref{uconv})
        \begin{align*}
            \frac{\rho_*}2\norm{\bu(t_n)}^2_{\bLds}\leq \frac12\int_\Omega\rho(\vphi(t_n))\normmm{\bu(t_n)}^2\dx\leq \frac{\rho^*}2\norm{\bu(t_n)}^2_{\bLds}\to0,\quad\text{ as }n\to\infty.
        \end{align*}
        The proof is concluded.
\end{proof}

   Using Lemma \ref{convaaa1}, upon noticing that the equation of the chemical potential $\mu$ \textit{does not} depend explicitly on the velocity $\bu$, we can follow \textit{verbatim} the same proof of Lemma \ref{twoparts}, and deduce that the very same lemma (which we do not rewrite for the sake of brevity) also holds under the assumptions of Theorem \ref{weakAGG}.

   Observe that properties \eqref{sepaglobal1}-\eqref{convu1} are a consequence of Lemma \ref{lemmauconv} and Lemma \ref{convaaa1}. Then, we are left to show the fact that the $\omega$-limit is a singleton. This can be done arguing as in the proof of Theorem \ref{uniqueeq}. Here we highlight the main differences. First, we notice that till inequality \eqref{pp} the same results also hold in this case, and we refer to the proof of Theorem \ref{uniqueeq} for the notation. The difference lies in the fact that here the energy inequality accounts for the velocity $\bu$ as well. Namely, recalling that $m(\cdot)\geq m_*>0$ and $\nu(\cdot)\geq \nu_*>0$, and setting $C_*=\min\{m_*,\nu_*\}>0$, we find from \eqref{energyineq2} that
					\begin{align*}
						&\nonumber C_*\left(\int_s^t\norm{\nabla \mu(\tau)}_{\bLd}^2\d\tau+\int_s^t\norm{D\bu(\tau)}^2_{\bLd}\d\tau\right)\\&\leq \int_s^t\int_\Omega \nu(\vphi(x,\tau))\normmm{D\bu(x,\tau)}^2\dx\d\tau+\int_s^t \int_\Omega m(\vphi(x,\tau))\normmm{\nabla\mu(x,\tau)}^2\dx\dtau
						\\&\leq E_{tot}(\bu(s),\vphi(s))-E_{tot}(\bu(t),\vphi(t)),
					\end{align*}
					for any $t>0$ and almost any $s\in[0,t]$, $s=0$ included.
					This gives
					\begin{align}
						& \nonumber\left(\int_s^t \norm{\nabla\mu(\tau)}_{\bLd}^2\dtau+\int_s^t\norm{D\bu(\tau)}^2_{\bLd}\d\tau\right)^{2(1-\vartheta)}
						\\&\leq \frac1{C_*^{2(1-\vartheta)}}(E_{tot}(\bu(s),\vphi(s))-E_{tot}(\bu(t),\vphi(t)))^{2(1-\vartheta)}.\label{cv1}
					\end{align}
					Recalling the proof of Lemma \ref{convaaa1}, we know that $E_{tot}(\bu(t),\vphi(t))\to E_\infty$ as $t\to \infty$. Then, passing to the limit as $t\to\infty$ in \eqref{cv1}, and then using Korn's inequality (with constant $C_K$), we get, for almost any $s\in(t_*,\infty)$,
					\begin{align}
						&\nonumber \left(\int_s^\infty \norm{\nabla\mu(\tau)}_{\Ld}^2\dtau+\int_s^\infty\norm{D\bu(\tau)}^2_{\bLd}\d\tau\right)^{2(1-\vartheta)}\\&
						\leq \frac1{C_*^{2(1-\vartheta)}}(E_{tot}(\vphi(s))-E_\infty)^{2(1-\vartheta)}\nonumber\\&\nonumber\leq \frac{2^{2(1-\vartheta})}{C_*^{2(1-\vartheta)}}\normmm{E(\vphi(s))-E_\infty}^{2(1-\vartheta)}(\chi_{A_M(T_S)}(s)+\chi_{(t_*,\infty)\setminus A_M(T_S)}(s))+2^{2(1-\vartheta)}\left(\frac{\rho^*}2\right)^{2(1-\vartheta)}\norm{\bu(s)}^{4(1-\vartheta)}_{\bLds}\\&
                        \leq \frac{2^{2(1-\vartheta})}{C_*^{2(1-\vartheta)}}\normmm{E(\vphi(s))-E_\infty}^{2(1-\vartheta)}(\chi_{A_M(T_S)}(s)+\chi_{(t_*,\infty)\setminus A_M(T_S)}(s))\nonumber\\&\quad +2^{2(1-\vartheta)}\left(\frac{\rho^*}2\right)^{2(1-\vartheta)}\frac{1}{C_*^{2(1-\vartheta)}}\left(\frac{2}{\rho_*}E_{tot}(\bu_0,\vphi_0)\right)^{1-2\vartheta}\norm{\bu(s)}^{2}_{\bLds}\nonumber\\&
                         \leq \frac{2^{2(1-\vartheta})}{C_*^{2(1-\vartheta)}}\normmm{E(\vphi(s))-E_\infty}^{2(1-\vartheta)}(\chi_{A_M(T_S)}(s)+\chi_{(t_*,\infty)\setminus A_M(T_S)}(s))\nonumber\\&\quad +2^{2(1-\vartheta)}\left(\frac{\rho^*}2\right)^{2(1-\vartheta)}
                         \frac{C_K^2}{C_*^{2(1-\vartheta)}}\left(\frac{2}{\rho_*}E_{tot}(\bu_0,\vphi_0)\right)^{1-2\vartheta}
                         \norm{D\bu(s)}^{2}_{\bLds}.\label{cv22b}
					\end{align}
                    Here we used the control $\sup_{t\geq0}\norm{\bu(t)}_{\bLds}^2\leq \frac{2}{\rho_*}E_{tot}(\bu_0,\vphi_0)$ which comes from the energy inequality \eqref{energyineq2}.
				Note that, by Lemma \ref{twoparts} (which also holds in this case), for almost any $s\in (t_*,\infty)\setminus A_M(T_S)$ it holds $\norm{\nabla \mu(s)}_{\bLd}\geq M$. Thus, recalling that $E(\vphi(t))\leq E_{tot}(\bu_0,\vphi_0)$ for any $t\geq0$, as well as $E_\infty\leq E_{tot}(\bu_0,\vphi_0) $, we get
					$$
				\normmm{E(\vphi(s))-E_\infty}^{2(1-\vartheta)}\chi_{(t_*,\infty)\setminus A_M(T_S)}(s)\leq (2E_{tot}(\bu_0,\vphi_0))^{2(1-\vartheta)}\frac{\norm{\nabla\mu(s)}_{\bLd}^2}{M^2}\chi_{(t_*,\infty)\setminus A_M(T_S)}(s),
				$$
				for almost any $s\in (t_*,\infty)\setminus A_M(T_S)$.
				Then, exploiting \eqref{pp}, we have
					$$
				\normmm{E(\vphi(s))-E_\infty}^{2(1-\vartheta)}\chi_{ A_M(T_S)}(s)\leq C\frac{\norm{\nabla\mu(s)}_{\bLd}^2}{M^2}\chi_{A_M(T_S)}(s),
				$$
				for any $s\in  A_M(T_S)$.
				As a consequence, we deduce from \eqref{cv22b} that
					\begin{align}
					& \left(\int_s^\infty \norm{\nabla\mu(\tau)}_{\bLd}^2\dtau+\int_s^\infty \norm{D\bu(\tau)}_{\bLd}^2\dtau\right)^{2(1-\vartheta)}
					\nonumber\\&\leq \frac{2^{2(1-\vartheta)}}{C_*^{2(1-\vartheta)}}\left(C+\frac{(2 E(\vphi_0))^{2(1-\vartheta)}}{M^2}\right)\norm{\nabla \mu (s)}^2_{\bLd}\nonumber\\&\nonumber\quad +2^{2(1-\vartheta)}\left(\frac{\rho^*}2\right)^{2(1-\vartheta)}\frac{C_K^2}{C_*^{2(1-\vartheta)}}\left(\frac{2}{\rho_*}E_{tot}(\bu_0,\vphi_0)\right)^{1-2\vartheta}\norm{D\bu(s)}^{2}_{\bLds}\nonumber\\&
                    \leq 2^{2(1-\vartheta)}\max\left\{\frac{1}{C_*^{2(1-\vartheta)}}\left(C+\frac{(2 E_{tot}(\bu_0,\vphi_0))^{2(1-\vartheta)}}{M^2}\right),\left(\frac{\rho^*}2\right)^{2(1-\vartheta)}\frac{C_K^2}{C_*^{2(1-\vartheta)}}\left(\frac{2}{\rho_*}E_{tot}(\bu_0,\vphi_0)\right)^{1-2\vartheta}\right\}\nonumber\\&\quad\times\left(\norm{\nabla\mu(s)}_{\bLd}^2+\norm{D\bu(s)}_{\bLd}^2\right)=:C_M \left(\norm{\nabla\mu(s)}_{\bLd}^2+\norm{D\bu(s)}_{\bLd}^2\right)
                    ,\label{cv2b}
				\end{align}
				for almost any  $s\in (t_*,\infty)$.
			    We now use Lemma \ref{Feireisl} once more with
$$
Z(\cdot)=\norm{\nabla\mu(\cdot)}_{\bLd}, \quad \alpha=2(1-\vartheta)\in(1,2), \quad \zeta=C_M>0,
$$
and $\mathcal M=(t_*,\infty)$. This entails that (see \eqref{fin})
					\begin{align}
						\nabla \mu\in L^1(t^*,\infty;L^2(\Omega)),\quad D\bu\in L^1(t^*,\infty;\bLd).\label{bA2}
					\end{align}
                    As a consequence, by Korn's inequality, recalling that $\normmm{\vphi}< 1$ almost everywhere in $\Omega\times(0,\infty)$, we see that
\begin{align*}
&\normmm{\int_\Omega \nabla\vphi\cdot \bu v\dx}=\normmm{\int_\Omega \vphi\bu\cdot   \nabla v\dx}\\
&\leq \norm{\vphi}_{L^\infty(\Omega)}\norm{\bu}_{\bLds}\normm{v}_{H^1(\Omega)}\leq \norm{\bu}_{\bLds}\normm{v}_{H^1(\Omega)}\\
&\leq C_K\norm{D\bu}_{\bLd}\normm{v}_{H^1(\Omega)},\quad \forall v\in H^1(\Omega).
\end{align*}
Therefore, recalling \eqref{bA2}, we get
\begin{align*}
    \nabla\vphi\cdot \bu\in L^1(t^*,\infty;H^1(\Omega)'),
\end{align*}
so that, by comparison in \eqref{CHeq}, we obtain $\partial_t\vphi\in L^1(t^*,\infty;\Hu')$. Hence, we have
$$
\vphi(t)=\vphi(t^*)+ \int_{t^*}^t\partial_t\vphi(\tau) \: \d\tau \to {\vphi_\infty}\quad \text{ in }\Hu', \text{ as } t \to \infty,
$$
for some ${\vphi_\infty}\in \Hu'$. As such, we have that $\vphi(t)$ converges in $\Hu'$ as $t\to\infty$ and, by uniqueness of the limit, recalling \eqref{uconv}, we conclude that $\omega(\bu,\vphi)$ is a singleton, i.e., $\omega(\bu,\vphi)=\{(\mathbf 0,\vphi_\infty)\}$. The convergence \eqref{equil1}  then follows from \eqref{convergence2}, as now $\omega(\bu,\vphi)$ is a singleton. The proof is then complete.

	\appendix
	\section{Some technical lemmas}
	\subsection{A Lemma on the integrability of functions}
	The following lemma, whose proof can be found in \cite[Lemma 7.1]{FS}, guarantees that a function is integrable if it satisfies a suitable integral inequality.		
			
\begin{lemma}
	Let $Z\geq 0$ be a measurable function on $(0,\infty)$ such that
	\begin{align*}
		Z\in L^2(0,\infty),\quad \int_0^\infty\normmm{Z(t)}^2\dt\leq Y,
	\end{align*}
	for some $Y>0$. If there exist $\alpha\in (1, 2)$, $\zeta>0$, and an open set $\mathcal M\subset (0,\infty)$ such that
	\begin{align}
		\left(\int_s^\infty Z^2(t)\dt\right)^{\alpha}\leq \zeta Z^2(s),\quad \text{ for a.a. }s\in \mathcal M,
	\end{align}
	then $Z\in L^1(\mathcal M)$, and there exists $C=C(Y,\alpha,\zeta)>0$, independent of $\mathcal M$, such that
	\begin{align}
		\int_{\mathcal M} Z(t)\dt \leq C.
		\label{fin}
	\end{align}
	\label{Feireisl}
\end{lemma}
\subsection{A lemma on geometric convergence of sequences}
 This lemma, which is a key tool in De Giorgi iteration argument, can be found, e.g., in \cite[Ch. I, Lemma 4.1]{DiBenedetto} (we refer to \cite[Lemma 3.8]{P} for a proof).
\begin{lemma}
	\label{conv}
	Let $\{y_n\}_{n\in\N\cup \{0\}}\subset \R^+$ satisfy the recursive inequalities
	\begin{align}
		y_{n+1}\leq Cb^ny_n^{1+\varepsilon},
		\label{ineq}\qquad \forall n\geq 0,
	\end{align}
	for some $C>0$, $b>1$ and $\varepsilon>0$. If
	\begin{align}
		\label{condition}
		y_0\leq \theta:= C^{-\frac{1}{\varepsilon}}b^{-\frac{1}{\varepsilon^2}},
	\end{align}
	then
	\begin{align}
		y_n\leq \theta b^{-\frac{n}{\varepsilon}},\qquad \forall n\geq 0,
		\label{yn}
	\end{align}
	and consequently $y_n\to 0$ for $n\to \infty$.
\end{lemma}

\medskip
  \textbf{Acknowledgments.}
   This research was funded in part by the Austrian Science Fund (FWF) \href{https://doi.org/10.55776/ESP552}{10.55776/ESP552}.
AP and MG are also members of Gruppo Nazionale per l’Analisi Matematica, la Probabilità e le loro Applicazioni (GNAMPA) of
Istituto Nazionale per l’Alta Matematica (INdAM). This research is part of the activities of “Dipartimento di Eccellenza 2023-2027” of
Politecnico di Milano (MG).
For open access purposes, the authors have applied a CC BY public copyright license to
any author accepted manuscript version arising from this submission.
\\


\begin{thebibliography}{30}


\bibitem{Abels}
{\sc H.~Abels}, {\em On a diffuse interface model for two-phase flows of viscous, incompressible fluids with matched densities},
Arch. Ration. Mech. Anal., 194 (2009), 463--506.



\bibitem{ADG}
{\sc H.~Abels, D.~Depner, and H.~Garcke}, {\em Existence of weak solutions for
  a diffuse interface model for two-phase flows of incompressible fluids with
  different densities}, J. Math. Fluid Mech., 15 (2013),  453--480.


\bibitem{AGG} {\sc H. Abels, H. Garcke and G. Gr\"un}, {\em Thermodynamically consistent, frame indifferent diffuse interface models for incompressible two-phase flows with different densities}, Math. Models Methods Appl. Sci. { 22} (2012), Paper No.~3, 1150013, 40 pp.

\bibitem{AGGG}
{\sc H.~Abels, H.~Garcke, and A.~Giorgini}, {\em Global regularity and
  asymptotic stabilization for the incompressible
  {N}avier--{S}tokes-{C}ahn--{H}illiard model with unmatched densities}, Math.
  Ann., 389 (2024), 1267--1321.

\bibitem{AGGmulti}
{\sc H.~Abels, H.~Garcke, and A.~Poiatti}, {\em Mathematical analysis of a
  diffuse interface model for multi-phase flows of incompressible viscous
  fluids with different densities}, J. Math. Fluid Mech., 26 (2024), Paper
  No. 29, 51 pp.


\bibitem{AGGP1}{\sc H. Abels, H. Garcke and A. Poiatti}, {\em Diffuse interface model for two-phase flows on evolving surfaces with different densities: global well-posedness}, Calc. Var. Partial Differential Equations { 64} (2025), no.~5, Paper No. 141, 41 pp.


\bibitem{AGGP2}{\sc H. Abels, H. Garcke and A. Poiatti}, {\em Diffuse interface model for two-phase flows on evolving surfaces with different densities: local well-posedness}, J. Evol. Equ., to appear (2026), 39 pp.


\bibitem{AGP1}
{\sc H.~Abels and A.~Poiatti}, {\em Weak solutions to a sharp interface model
  for a two-phase flow of incompressible viscous fluids with different
  densities}, arXiv:2505.06423v1 [math.AP] (2025), 56 pp.

\bibitem{AW}
{\sc H.~Abels and M.~Wilke}, {\em Convergence to equilibrium for the
  {C}ahn-{H}illiard equation with a logarithmic free energy}, Nonlinear Anal.,
  67 (2007), 3176--3193.

\bibitem{Alberti}
{\sc S.~Alberti}, {\em Phase separation in biology}, Curr. Biol., 27 (2017), R1097--R1102.

\bibitem{BB}
{\sc J.W.~Barrett and J.F.~Blowey}, {\em Finite element approximation of the
  {C}ahn-{H}illiard equation with concentration dependent mobility}, Math.
  Comp., 68 (1999), 487--517.

\bibitem{CEGP}
{\sc D.~Caetano, C.M.~Elliott, M.~Grasselli, A.~Poiatti}, {\em Regularization and separation for evolving surface Cahn-Hilliard equations},
SIAM J. Math. Anal., 55 (2023), 6625--6675.

\bibitem{Cahn}
{\sc J.W.~Cahn}, {\em On spinodal decomposition}, Acta Metallurgica, 9 (1961), 795--801.

\bibitem{CH1}
{\sc J.W.~Cahn, J.E.~Hilliard}, {\em Free energy of a nonuniform system. I. Interfacial free energy},
J. Chem. Phys., 28 (1958) 258--267.


\bibitem{CH2}
{\sc J.W.~Cahn, J.E.~Hilliard},
{\em Spinodal decomposition: a reprise}, Acta Metallurgica, 19 (1971), 151--161.


\bibitem{CGGG}
{\sc M.~Conti, P.~Galimberti, S.~Gatti, and A.~Giorgini}, {\em New results for
  the {C}ahn-{H}illiard equation with non-degenerate mobility: well-posedness
  and longtime behavior}, Calc. Var. Partial Differential Equations, 64 (2025),
  Paper No. 87, 32 pp.

\bibitem{DiBenedetto}
{\sc E.~DiBenedetto}, {\em Partial Differential Equations}, Cornerstones,
  Birkh\"{a}user Boston, MA, 2nd edition~ed., 2009.

\bibitem{DiPG}
{\sc A.~Di Primio, M. Grasselli}, {\em Analysis of a diffuse interface model for two-phase magnetohydrodynamic flows},
Discrete Contin. Dyn. Syst. Ser. S, 16 (2023), 3473--3534.

\bibitem{Dolgin}
{\sc E.~Dolgin}, {\em How phase separation is revolutionizing biology}, Nat., 626 (2024), 1152--1154.

\bibitem{Elliott}
{\sc C.M. Elliott}, {\em The Cahn--Hilliard model for the kinetics of phase separation}, in Mathematical models for phase change problems
(\'{O}bidos, 1988), 35--73, Internat. Ser. Numer. Math. 88, Birkh\"{a}user, Basel, 1989.

\bibitem{EG}
{\sc C.M. Elliott, H. Garcke}, {\em On the Cahn--Hilliard equation with degenerate mobility},
SIAM J. Math. Anal., 27 (1996), 404--423.

\bibitem{FS}
{\sc E.~Feireisl and F.~Simondon}, {\em Convergence for semilinear degenerate
  parabolic equations in several space dimensions}, J. Dynam. Differential
  Equations, 12 (2000), 647--673.

\bibitem{FHLS}
{\sc J.~Fischer, S.~Hensel, T.~Laux, and T.~Simon}, {\em A weak-strong
  uniqueness principle for the {M}ullins-{S}ekerka equation}, arXiv:2404.02682v1 [math.AP]
   (2024), 48 pp.
   
   \bibitem{FrigeriGrasselli}
   {\sc S.~Frigeri and M.~Grasselli}, {\em Nonlocal Cahn--Hilliard-Navier--Stokes
   	systems with singular potentials,} Dyn. Partial Differ. Equ., 9 (2012), 273--304.

\bibitem{GGPS}
{\sc C.G.~Gal, M.~Grasselli, A.~Poiatti, and J.L. Shomberg}, {\em
  Multi-component {C}ahn-{H}illiard systems with singular potentials:
  theoretical results}, Appl. Math. Optim., 88 (2023), Paper No. 73, 46 pp.




\bibitem{GalP}
{\sc C.G.~Gal and A.~Poiatti}, {\em Unified framework for the separation
  property in binary phase-segregation processes with singular entropy
  densities}, European J. Appl. Math., 36 (2025), 40--67.

\bibitem{Gio}
{\sc A.~Giorgini},
{\em Well-posedness of a diffuse interface model for Hele-Shaw flow},
J. Math. Fluid Mech., 22 (2020), Paper No. 5, 36 pp.

\bibitem{GPCAC}
{\sc M.~Grasselli and A.~Poiatti}, {\em Multi-component conserved
  {A}llen--{C}ahn equations}, Interfaces Free Bound., 26 (2024), 489--541.

\bibitem{Haraux}
{\sc A.~Haraux}, {\em Syst\`{e}mes Dynamiques Dissipatifs et Applications},
Rech. Math. Appl. 17, Masson, Paris, 1991.

\bibitem{HeWu}
{\sc J.~He, H.~Wu}, {\em On a Navier--Stokes-Cahn--Hilliard system for viscous incompressible two-phase flows with chemotaxis, active transport and reaction},
Math. Ann., 389 (2024), 2193--2257.

\bibitem{HKP}
{\sc C.~Hurm, P.~Knopf, and A.~Poiatti}, {\em Nonlocal-to-local convergence
  rates for strong solutions to a {N}avier-{S}tokes-{C}ahn--{H}illiard system
  with singular potential}, Commun. Partial Differential Equations, 49
  (2024), 832–871.

\bibitem{Kenmochi}
{\sc N.~Kenmochi, M.~Niezg\'odka, and I.~Paw\l~ow}, {\em Subdifferential
  operator approach to the {C}ahn--{H}illiard equation with constraint}, J.
  Differential Equations, 117 (1995), 320--356.

\bibitem{Miranville}
{\sc A. Miranville}, {\em The {C}ahn--{H}illiard Equation: Recent Advances and Applications}, CBMS-NSF Regional Conf.
Ser. in Appl. Math. 95. Society for Industrial and Applied Mathematics, Philadelphia, PA, 2019.

\bibitem{P}
{\sc A.~Poiatti}, {\em The 3{D} strict separation property for the nonlocal
  {C}ahn--{H}illiard equation with singular potential}, Anal. PDE, 18 (2025),
  109--139.

\bibitem{PolSim}
{\sc P.~Pol\'{a}\v{c}ik, F.~Simondon},
{\em Nonconvergent bounded solutions of semilinear heat equations on arbitrary domains},
J. Differential Equations, 186 (2002), 586--610.

\bibitem{RybHof}
{\sc P.~Rybka, K.-H.~Hoffmann}, {\em Convergence of solutions to {C}ahn-{H}illiard equation},
Comm. Partial Differential Equations, 24 (1999), 1055--1077.

\bibitem{Schimperna}
{\sc G.~Schimperna}, {\em Global attractors for {C}ahn-{H}illiard equations
  with nonconstant mobility}, Nonlinearity, 20 (2007), 2365--2387.

\bibitem{Wu}
{\sc H.~Wu}, {\em A review on the {C}ahn-{H}illiard equation: classical results and recent advances in dynamic boundary conditions},
Electron. Res. Arch., 30 (2022), 2788--2832.

\end{thebibliography}


\end{document}